\documentclass[a4paper, 11pt]{article}
\usepackage{amsmath} \usepackage{amsfonts} \usepackage{amssymb}\usepackage{latexsym}
\usepackage[english]{babel}
\usepackage{amscd}
\usepackage[dvips,final]{graphics}
\usepackage{locengli,math}
\usepackage{graphicx,pst-plot,psfig,pstricks,pst-node,pst-text,pst-tree,pst-char,pst-coil}
\begin{document}
\begin{center}
\textbf{\LARGE{\textsf{Ennea-algebras}}}
\footnote{
{\it{2000 Mathematics Subject Classification: 17A30, 18D50.}}
{\it{Key words and phrases: Nested dendriform trialgebras, $t$-ennea-algebras, quadri-algebras, $t$-Baxter operators, $t$-infinitesimal bialgebras, left pre-Lie algebras, operads, connected Hopf algebras, formal deformations.}} 
}
\vskip1cm
\parbox[t]{14cm}{\large{
Philippe {\sc Leroux}}\\
\vskip4mm
{\footnotesize
\baselineskip=5mm
Institut de Recherche
Math\'ematique, Universit\'e de Rennes I and CNRS UMR 6625\\
Campus de Beaulieu, 35042 Rennes Cedex, France, pleroux@univ-rennes1.fr}}
\end{center}

\vskip1cm
{\small
\vskip1cm
\baselineskip=5mm
\noindent
{\bf Abstract:}
A generalisation of a recent work of M. Aguiar and J.-L. Loday on quadri-algebras called $t$-ennea-algebras constructed over dendriform trialgebras is proposed. Such algebras allow the construction of nested dendriform trialgebras and are related to pre-Lie algebras, $t$-infinitesimal bialgebras and $t$-Baxter operators. 
We also show that the augmented free $t$-ennea-algebra has a structure of connected Hopf algebra. In the last part, we use Baxter operators to produce formal deformations of dendriform algebras, quadri-algebras and ennea-algebras.
Examples and relations of this work to combinatorial objects are given.
\begin{scriptsize}
\begin{center}
\tableofcontents
\end{center}
\end{scriptsize}
\section{Introduction}
In this paper, $k$ denotes a field of characteristic zero. 

\noindent
\textbf{Notation:}
Let $(X, \ \bullet)$ be a $k$-algebra and $(\bullet_i)_{1 \leq i \leq N}: X^{\otimes 2} \xrightarrow{} X$ be a family of binary operations on $X$. 
The notation $\bullet \longrightarrow \sum_i \bullet_i$ will mean $x \bullet y =  \sum_i x \bullet_i y$, for all $x,y \in X$. We say that
the operation $\bullet$ {\it{splits}} into the $N$ operations $\bullet_1, \ldots, \bullet_N,$ or that the operation $\bullet$ is a {\it{cluster of $N$ (binary) operations}}. If $\bullet$ is associative then such a cluster will be said {\it{associative}}.

Motivated by $K$-theory, J.-L. Loday first introduced  a ``non-commutative version'' of Lie algebras called Leibniz algebras \cite{LodayLeib}. Such algebras are described by a bracket $[ -,z]$ verifying the Leibniz identity:
$$[[x,y ] ,z ] = [[ x,z ] ,y ] + [ x,[ y  ,z ]].$$
When the bracket is skew-symmetric, the Leibniz identity becomes the
Jacobi identity and Leibniz algebras turn out to be Lie algebras. A way to construct such  Leibniz algebras is to start with  associative dialgebras, i.e., $k$-vector spaces $D$ equipped with two associative products,
$\vdash$ and $\dashv$, and verifying some conditions \cite{Loday}. The operad $Dias$ associated with associative dialgebras is then Koszul dual to the operad $DiDend$ associated with dendriform dialgebras \cite{Loday}.
A {\it{dendriform dialgebra}} is a $k$-vector space $E$ equipped with two binary operations:
$\prec, \ \succ: E^{\otimes 2} \xrightarrow{} E$, satisfying the following relations for all $x,y \in E$:
$$(x \prec y )\prec z = x \prec(y \star z), \ \ \
(x \succ y )\prec z = x \succ(y \prec z), \ \ \
(x \star y )\succ z = x \succ(y \succ z), \ \ \ $$
where, by definition, $x \star y :=x  \prec y +x \succ y$, for all $x,y \in E$. The $k$-vector space $(E,\star)$ is then an associative algebra such that $\star \longrightarrow \ \prec + \succ$. Observe that these axioms are globally invariant under the transformation $x \prec^{op} y := y \succ x$, $x \succ^{op} y := y \prec x$. A dendriform dialgebra is said to be commutative if $x \prec y := y \succ x$. Such algebras are also called Zinbiel algebras \cite{Loday}.
Similarly, to propose a ``non-commutative version'' of Poisson algebras, J.-L. Loday and M. Ronco \cite{LodayRonco} introduced the notion of associative trialgebras. It turns out that $Trias$, the operad associated with this type of algebras, is Koszul dual to $TriDend$, the operad associated with
dendriform trialgebras.
\begin{defi}{[Dendriform trialgebra]}
\label{tridend}
A {\it{dendriform trialgebra}} is a $k$-vector space $T$ equipped with three binary operations:
$\prec, \ \succ, \ \circ: T^{\otimes 2} \xrightarrow{} T$, satisfying the following relations for all $x,y \in T$:
$$(x \prec y )\prec z = x \prec(y \star z), \ 
(x \succ y )\prec z = x \succ(y \prec z), \
(x \star y )\succ z = x \succ(y \succ z), $$
$$(x \succ y )\circ z = x \succ(y \circ z), \ 
(x  \prec y )\circ z = x \circ(y \succ z), \ 
(x \circ y )\prec z = x \circ(y \prec z), \
(x \circ y ) \circ z = x \circ(y \circ z),$$
where by definition $x \star y :=x  \prec y +x \succ y +x \circ y$, for all $x,y \in T$. The $k$-vector space $(T,\star)$ is then an associative algebra such that $\star \longrightarrow \ \prec + \succ + \ \circ$.
\end{defi}
Observe that these axioms are globally invariant under the transformation $x \prec^{op} y := y \succ x$, $x \succ^{op} y := y \prec x$ and $x \circ^{op} y := y \circ x$. A dendriform trialgebra is said to be {\it{commutative}} if $x \prec y := y \succ x$ and $x \circ y := y \circ x$.
\NB
As the product $\circ$ is associative so is the product $t \circ$ defined by $x (t \circ) y := t x  \circ y$, where $t \in k$ and $x,y \in T$. Therefore, if 
$(T, \ \prec, \ \succ, \ \circ)$ is a dendriform trialgebra so is $(T, \ \prec, \ \succ, \ t\circ)$. For $t=0$, the dendriform trialgebra $(T, \ \prec, \ \succ, \ 0\circ)$ is a dendriform dialgebra. 

In \cite{Aguiar}, M. Aguiar showed important relations between Baxter operators and infinitesimal bialgebras, which gave birth to  {\it{quadri-algebras}} in \cite{AguiarLoday}. A {\it{quadri-algebra}} is a $k$-vector space $Q$ equipped with four binary operations $\nwarrow, \ \nearrow, \ \swarrow, \ \searrow \ : \ Q^{\otimes 2} \xrightarrow{} Q$ and obeying compatibility axioms. One of the interests of such algebras is to produce an associative algebra $(Q, \star)$ such that $\star$ is the following associative cluster $\star  \xrightarrow{} \ \nwarrow + \nearrow + \swarrow + \searrow$. The second and main interest is to produce two dendriform dialgebras by combining these four operations, the two dendriform dialgebras $(Q, \prec, \succ)$ and $(Q, \wedge, \vee)$ being such that $\star \xrightarrow{} \ \prec + \succ$ and $\star \xrightarrow{} \  \wedge + \vee$.
Quadri-algebras can be constructed via Baxter operators and infinitesimal bialgebras.

What has been done for dendriform dialgebras and quadri-algebras can be also realised for dendriform trialgebras and what are called $t$-ennea-algebras \footnote{In Greek, ennea means nine.}, $t \in k$. Let $(T, \star)$ be an associative algebra, such that $\star$ is an associative cluster, $\star \xrightarrow{} \ \prec + \succ + \ \circ$, of three laws $\prec$, $\succ$ and $\circ$ obeying the axioms of Definition \ref{tridend}. We construct another type of algebra, called {\it{$t$-ennea-algebra}}, via $t$-Baxter operators, $t \in k$.
Such an algebra, constructed over the dendriform trialgebra $(T, (\star), \ \prec, \ \succ, \ \circ)$ is equipped with six binary operations: 
$\nwarrow, \ \nearrow, \ \uparrow, \ \searrow, \ \swarrow, \ \downarrow \ : \ T^{\otimes 2} \xrightarrow{} T$, verifying compatibility axioms. Combining these nine operations in a suitable way gives an associative algebra whose associative product $\bar{\star}$ splits into three laws in two ways: $\bar{\star} \xrightarrow{} \ \triangleright + \triangleleft + \ \bar{\circ}$ and $\bar{\star} \xrightarrow{} \ \vee + \wedge + \ t \star$, the algebraic structures $(T, \ \triangleleft , \ \triangleright, \ \bar{\circ} )$ and $(T, \ \wedge, \ \vee, \ t\star)$ being both dendriform trialgebras. The novelty is that the associative algebras $(T, \ \bar{\circ})$ and $(T, \ t \star)$ are also both dendriform trialgebras, i.e., their associative laws are clusters of three operations verifying the axioms of Definition \ref{tridend}.

A dendriform trialgebra $(T, \ \prec, \ \succ, \ \ \circ)$ is said to be \textit{nested} if the associative product of the associative algebra $(T, \ \circ)$ is itself a cluster of three laws verifying  dendriform trialgebra axioms. A way to produce such nested 
dendriform trialgebras is to construct $t$-ennea-algebras, themselves constructed via $t$-infinitesimal bialgebras and $t$-Baxter operators. For the case $t=0$, the framework of Aguiar and Loday is recovered. 

Section~2 gives the definition of $t$-Baxter operators as well as the axioms of  $t$-ennea-algebras.
These new structures give rise to (left) pre-Lie algebras. In particular, we show that a $t$-ennea-algebra gives two pre-Lie algebras whose Lie algebra structures coincide.
Section~3 relates $t$-ennea-algebras to $t$-infinitesimal bialgebras. Section~4 relates $t$-infinitesimal bialgebras to (left) pre-Lie algebras. Examples from weighted graph theory are given in Section~5. In Section~6, we say a word on the free $t$-ennea-algebra on a $k$-vector space $V$ and combinatorial objects related to it. In particular, we show that the augmented free $t$-ennea-algebra has a connected Hopf algebra structure. This work ends with formal deformations in Section~7. We show the usefulness of Baxter operators and 2-hypercubic coalgebras \cite{codialg1} to produce formal deformations of dendriform algebras, quadri-algebras and ennea-algebras. Explicit examples are given.
\section{On $t$-Baxter operators}
\subsection{$t$-Baxter operators and $t$-ennea-algebras}
\begin{defi}{[$t$-Baxter operator]}
Let $t \in k$.
Let $(A, \ \cdot)$ be an associative algebra. A {\it{$t$-Baxter operator}} is a linear map $\beta: A \xrightarrow{} A$ verifying:
$$\beta(x) \cdot  \beta(y) = \beta(x \cdot \beta(y) + \beta(x)\cdot y + tx \cdot y).  $$
\end{defi}
\Rk
For $t=0$, this map is called a Baxter operator. It
appears originally in a work of G. Baxter \cite{Baxter} and the importance of such a map was stressed by G.-C. Rota in \cite{Rota1}.
If $t=-1$, such a map plays the r\^ole of the renormalisation map in quantum field theory in the work of A. Connes and D. Kreimer \cite{Kreimer, CKreimer}. Let us give examples.
\begin{prop} \textsf{[Change of scale]} \ \
\label{Scale}
Let $t,t' \in k$.
Let $(A, \ \cdot)$ be an associative algebra and
$\beta: A \xrightarrow{} A$, a $t$-Baxter operator. 
Via the transformation $\beta \mapsto \beta' := t'\beta$, the linear map $\beta'$ becomes a $tt'$-Baxter operator.
\end{prop}
\Proof
Straightforward.
\eproof
\begin{prop}
Denote by $M_n(k)$ the algebra of  $n$ by $n$ matrices and fix $t \in k$.
Let $T$ be the set of upper triangular matrices $T := \{X \in M_n(k); \ X = \begin{pmatrix}
 \cdot & *\\
0 & \cdot
\end{pmatrix} \}$. Then the linear map $\beta: T \xrightarrow{} T$ such that
$(\beta(X))_{ii} := t\sum_j X_{ij}$
and $(\beta(X))_{ij}=0$ for $i \not=j$,
is a $-t$-Baxter operator.
Similarly, the linear map $\beta_1: T \xrightarrow{} T$ such that
$(\beta_1(X))_{jj} := t\sum_i X_{ij}$
and $(\beta(X))_{ji}=0$ for $i \not=j$, is a $-t$-Baxter operator.
\end{prop}
\Proof
Straightforward.
\eproof

\noindent
Before going ahead, let us observe that the concept of $t$-Baxter operator has a straightforward dualisation. Let $(C, \ \Delta)$ be
a coassociative coalgebra. A {\it{$t$-coBaxter operator}}
is a linear map $\psi: C \xrightarrow{} C$ such that:
$$ (\psi \otimes \psi) \Delta = t\Delta\psi +  (id \otimes \psi)\Delta\psi + (\psi \otimes id )\Delta\psi. $$
For instance, define the coassociative coalgebra $C$
spanned by a basis $(X_{ij})_{1 \leq i \leq j \leq n}$
as a $k$-vector space and equipped with a coassociative
coproduct $\Delta$ defined by $\Delta X_{ij} := \sum_k X_{ik} \otimes X_{kj}$, for all $1 \leq i,j \leq n$ and a counit $\epsilon$ defined by $\epsilon(X_{ij}):= 0$, if $i \not= j$ and equals to 1 otherwise. Define the linear
map $\psi: C \xrightarrow{} C$ such that $\psi(X_{ij})=0$, if $i \not=j$ and $\psi(X_{ii})= t \sum_j X_{ij}$ otherwise,
then $\psi$ is a $(-t)$-coBaxter operator. In the sequel, we will focus on the algebra side, the dualisation being straightforward.

\noindent
The following proposition also appears in \cite{KEF}.
\begin{prop}
Let $(A, \ \cdot)$ be an associative algebra and $\beta: A \xrightarrow{} A$ be a $t$-Baxter operator. The following operations:
$$x \prec_{\beta} y := x \cdot \beta(y), \ x \succ_{\beta} y:= \beta(x) \cdot y, \ x \circ_{\beta} y = tx \cdot y,$$ 
determine a dendriform trialgebra structure on $A$. The map $\beta: (A, \star_{\beta}) \xrightarrow{} (A, \cdot)$, where $\star_{\beta}$
is defined such that $x \star_{\beta}  y:= x \prec_{\beta} y + x \succ_{\beta} y + x \circ_{\beta} y$, for all $x,y \in A$, is a morphism of associative algebras.
\end{prop}
\Proof
The proof is completed by checking axioms of Definition \ref{tridend}.
\eproof

Recall that a binary operation, associative or not, {\it{splits}} if
it can be decomposed into several binary operations, associative or not, like the associative product of dendriform trialgebras or dendriform dialgebras. Such a splitting can be obtained from the action of a $t$-Baxter operator, if it exists, on the associative algebra $(A, \ \cdot)$. This operator is an associative algebra morphism $(A, \ \star) \xrightarrow{}(A, \ \cdot)$ and the product $\star$ is an associative cluster of three operations obeying axioms of Definition \ref{tridend}. 

\noindent
From now on, we will always suppose $t \not= 0$, the case $t=0$ being studied in \cite{AguiarLoday}.
We now give a process, inspired from \cite{AguiarLoday}, to split the three products
of a dendriform trialgebra obtained, for instance, by an action of
a $t$-Baxter operator on the associative algebra $(A, \cdot)$.
\begin{defi}{[$t$-Ennea-algebra over a dendriform trialgebra]} 
\label{ennea}
Fix $t \in k$. A {\it{ $t$-ennea-algebra }} is a $k$-vector space $T$ equipped with nine binary operations
$ \prec, \ \succ, \ \circ, \ \uparrow, \ \downarrow, \ \searrow, \ \nearrow, \ \swarrow, \ \nwarrow: T \xrightarrow{} T^{\otimes 2}$ satisfying 49 relations.
To ease notation, consider the following sums defined for all $x,y \in T$ by,
\begin{eqnarray*}
x \triangleright y &:=& x \nearrow y + x \searrow y + t x \succ y, \\
x  \triangleleft y &:=& x \nwarrow y + x \swarrow y + t x \prec y, \\
x \bar{\circ} y &:=& x \uparrow y + x \downarrow y +  t x \circ y , \\
x \vee y &:=& x \searrow y + x \swarrow y +  x \downarrow y, \\
x \wedge y &:=& x \nearrow y + x \nwarrow y +  x \uparrow y, \\
x \star y &:=& x \succ y + x \prec y + x \circ y, 
\end{eqnarray*}
and $$ x \bar{\star} y := x \triangleright y + x  \triangleleft y + x \bar{\circ} y = x \wedge y + x \vee y + tx \star y.$$
The  
compatibility axioms are now presented in the following $7 \times 7$ matrix. 
\begin{footnotesize}
\begin{enumerate}
\item {
$$ 
(x \nwarrow y)\nwarrow z = x \nwarrow (y \bar{\star} z), \ \ \ \
(x \nearrow y)\nwarrow z = x \nearrow (y \triangleleft z), \ \ \ \
(x \wedge y)\nearrow z = x \nearrow (y \triangleright z), \ \ \ \ $$
$$
(x \nearrow y)\uparrow z = x \nearrow (y \bar{\circ} z), \ \ \ \ 
(x \nwarrow y)\uparrow z = x \uparrow (y \triangleright z), \ \ \ \
(x \uparrow y)\nwarrow z = x \uparrow (y \triangleleft z), \ \ \ \ 
(x \uparrow y)\uparrow z = x \uparrow (y \bar{\circ} z). \ \ \ \ 
$$ \\}
\item {
$$ 
(x \swarrow y)\nwarrow z = x \swarrow (y \wedge z), \ \ \ \
(x \searrow y)\nwarrow z = x \searrow (y \nwarrow z), \ \ \ \
(x \vee y)\nearrow z = x \searrow (y \nearrow z), \ \ \ \ $$
$$
(x \searrow y)\uparrow z = x \searrow (y \uparrow z), \ \ \ \ 
(x \swarrow y)\uparrow z = x \downarrow (y \nearrow z), \ \ \ \
(x \downarrow y)\nwarrow z = x \downarrow (y \nwarrow z), \ \ \ \ 
(x \downarrow y)\uparrow z = x \downarrow (y \uparrow z). \ \ \ \ 
$$ \\}
\item {
$$ 
(x \triangleleft y)\swarrow z = x \swarrow (y \vee z), \ \ \ \
(x \triangleright y)\swarrow z = x \searrow (y \swarrow z), \ \ \ \
(x \bar{\star} y)\searrow z = x \searrow (y \searrow z), \ \ \ \ $$
$$
(x \triangleright y)\downarrow z = x \searrow (y \downarrow z), \ \ \ \ 
(x \triangleleft y)\downarrow z = x \downarrow (y \searrow z), \ \ \ \
(x \bar{\circ} y)\swarrow z = x \downarrow (y \swarrow z), \ \ \ \ 
(x \bar{\circ} y)\downarrow z = x \downarrow (y \downarrow z).\ \ \ \ 
$$ \\}
\item {
$$ 
t(x \swarrow y)\prec z = tx \swarrow (y \star z), \ \ \ \ 
t(x \searrow y)\prec z = tx \searrow (y \prec z), \ \ \ \
t(x \vee y)\succ z = tx \searrow (y \succ z), \ \ \ \ $$
$$
t(x \searrow y)\circ z = tx \searrow (y \circ z), \ \ \ \
t(x \swarrow y)\circ z = tx \downarrow (y \succ z), \ \ \ \
t(x \downarrow y)\prec z = tx \downarrow (y \prec z), \ \ \ \ 
t(x \downarrow y) \circ z = tx \downarrow (y \circ z). \ \ \ \ 
$$ \\}
\item {
$$ 
t(x \nwarrow y)\prec z = tx \prec (y \vee z), \ \ \ \ 
t(x \nearrow y)\prec z = tx \succ (y \swarrow  z), \ \ \ \
t(x \wedge y)\succ z = tx \succ (y \searrow z), \ \ \ \ $$
$$
t(x \nearrow y)\circ z = tx \succ (y \downarrow z), \ \ \ \
t(x \nwarrow y)\circ z = tx \circ (y \searrow z), \ \ \ \
t(x \uparrow y)\prec z = tx \circ (y \swarrow z), \ \ \ \ 
t(x \uparrow y) \circ z = tx \circ (y \downarrow z).
$$ \\}
\item {
$$ 
t(x \prec y)\nwarrow z = tx \prec (y \wedge z), \ \ \ \ 
t(x \succ y)\nwarrow z = tx \succ (y \nwarrow  z), \ \ \ \
t(x \star y)\nearrow z = tx \succ (y \nearrow z), \ \ \ \ $$
$$
t(x \succ y)\uparrow z = tx \succ (y \uparrow z), \ \ \ \
t(x \prec y)\uparrow z = tx \circ (y \nearrow z), \ \ \ \
t(x \circ y)\nwarrow z = tx \circ (y \nwarrow z), \ \ \ \ 
t(x \circ y) \uparrow z = tx \circ (y \uparrow z), 
$$}
\item{
$$t^2(x \prec y)\prec z := t^2x \prec (y \star z), \ \ \ \
t^2(x \succ y)\prec z := t^2x \succ (y \prec z), \ \ \ \ 
t^2(x \star y)\succ z := t^2x \succ (y \succ z), \ \ \ \ $$
$$
t^2(x \succ y)\circ z := t^2x \succ (y \circ z), \ \ \ \
t^2(x \prec y)\circ z := t^2x \circ (y \succ z), \ \ \ \
t^2(x \circ y)\prec z := t^2x \circ (y \prec z), \ \ \ \
t^2(x \circ y)\circ z := t^2x \circ (y \circ z), \ \ \ \
$$}
\end{enumerate}
\end{footnotesize}
which give 49 relations. 
\Rk
Observe the last seven ones implies that $(T,\ \prec,\ \succ, \ \circ)$ is a dendriform trialgebra.
If the dendriform trialgebra $(T,\ \prec,\ \succ, \ \circ)$ can be ``augmented'' in this way, we shall say that the dendriform trialgebra
$T$ admits a {\it{ $t$-ennea-algebra structure}} or that $(T, \ \prec, \ \succ, \ \circ, \ \uparrow, \ \downarrow, \ \searrow, \ \nearrow, \ \swarrow, \ \nwarrow)$ is a {\it{$t$-ennea-algebra over the dendriform trialgebra}} $(T, \ \prec, \ \succ, \ \circ)$.
\Rk \textbf{[Transpose and opposite of a $t$-Ennea-algebra]}

\noindent
Observe also the existence of some symmetries.
The matrix remains globally invariant by interchanging $\swarrow$ and $\nearrow$; $t\succ$ and $\downarrow$; $\uparrow$ and $t\prec$; $\bar{\circ}$ and $t\star$. Indeed the operations
$\triangleright$ and $\vee$; $\triangleleft$ and $\wedge$, interchange. This new $t$-ennea-algebra is called the transpose of the previous $t$-ennea-algebra, as in \cite{AguiarLoday}. More generally, if $(T, \ \prec_1, \ \succ_1, \ \circ_1)$ and $(T, \ \prec_2, \ \succ_2, \ \circ_2)$ are two dendriform trialgebras admitting a $t$-ennea-algebra structures, we will say that $(T, \ \prec_2, \ \succ_2, \ \circ_2,\ \uparrow_2, \ \downarrow_2, \ \searrow_2, \ \nearrow_2, \ \swarrow_2, \ \nwarrow_2 )$ is the transpose of 
$(T, \ \prec_1, \ \succ_1, \ \circ_1, \ \uparrow_1, \ \downarrow_1, \ \searrow_1, \ \nearrow_1, \ \swarrow_1, \ \nwarrow_1 )$ if
$\swarrow_1=\nearrow_2$; $\nearrow_1=\swarrow_2$; $t\succ_1 = \downarrow_2$;
$t\succ_2 = \downarrow_1$; 
$\uparrow_1=t\prec_2$; 
$\uparrow_2=t\prec_1$;
$\bar{\circ}_1=t\star_2$; 
$\bar{\circ}_2=t\star_1$
and $\circ_1 =\circ_2$; $\nwarrow_1=\nwarrow_2$; $\searrow_1=\searrow_2$.
Similarly, we define the opposite of the $t$-ennea-algebra over a dendriform trialgebra
as follows:
$$x \searrow^{op} y = y \nwarrow x, \ \ x \nearrow^{op} y = y \swarrow x, \ \ x \nwarrow^{op} y = y \searrow x, \ \
x \swarrow^{op} y = y \nearrow x, \ \ x \uparrow^{op} y = y \downarrow x,$$
$$x \downarrow^{op} y = y \uparrow x, \ \ x \succ^{op} y =  y \prec x, \ \ x \circ^{op} y =y \circ x, \ \ x \prec^{op} y =  y \succ x .$$
Therefore, $x \triangleleft^{op} y = y \triangleright x, \ \ x \triangleright^{op} y = y \triangleleft x, \ \ x \vee^{op} y = y \wedge x, \ \ x \wedge^{op} y = y \vee x$.
A $t$-ennea-algebra is said to be {\it{commutative}} when it coincides with its opposite. Indeed, observe that for any $x,y \in T$, 
$x \bar{\star} y := y \bar{\star} x$.
\Rk
For $t=0$, the binary operations $t\prec$, $t \succ$, $t\circ$,  $\downarrow$ and $\uparrow$ vanish. The definition of quadri-algebras made out with four binary operations obeying 9 relations is then recovered.
\end{defi}
\begin{theo}
\label{enneathm}
Let $T$ be a
$t$-ennea-algebra. By using notation of Definition \ref{ennea}, 
the seven column sums of the matrix yield a dendriform trialgebra $(T_h,\ \triangleleft, \ \triangleright,\ \bar{\circ})$ called the horizontal dendriform trialgebra associated with $T$.
Similarly, the seven row sums of the matrix yield a dendriform trialgebra $(T_v,\ \wedge,\ \vee, \ t\star)$ called \footnote{Here, the law $t\star$ reads $x (t\star) y := t x \star y$, for all $x,y \in T$.} the vertical dendriform trialgebra associated with $T$.
\end{theo}
\Proof
Axioms of dendriform trialgebras are easily verified, for instance $(x \wedge y) \wedge z = x \wedge (y \ \bar{\star} \ z)$ and so forth.
\eproof

\noindent
Here are represented relations between categories:
\begin{center}
$
\begin{array}{ccc}
\textsf{Ennea.} & \longleftarrow & \textsf{Quadri.} \\ 
\downarrow & \ &  \downarrow  \\
\textsf{TriDend.} & \longleftarrow  & \textsf{DiDend.} \\
\downarrow & & \\
\textsf{As.} & & 
\end{array} 
$
\end{center}
And,
\begin{center}

$
\begin{array}{ccc}
\textsf{Ennea.} & \stackrel{F_v}{\longrightarrow} & \textsf{TriDend.} \\ 
 & & \\ 
 F_h \downarrow & \searrow f     & \downarrow F   \\
& & \\ 
\textsf{TriDend.} & \stackrel{F}{\longrightarrow}  & \textsf{As.}
\end{array} 
$
\end{center}
This diagram is commutative, i.e., $FF_v = f= FF_h$, where $F_v$ (resp. $F_h$) is the functor giving the vertical structure (resp. the horizontal structutre) of a ennea-algebra.
We now give a way, inspired from \cite{AguiarLoday}, to construct such $t$-ennea-algebras
over dendriform trialgebras. First of all, let us define the notion of a {\it{$t$-Baxter operator over
a dendriform trialgebra}}.
\begin{defi}{}
Let $(T, \ \star)$ be an associative algebra whose product $\star$ is an associative cluster turning $T$ into a dendriform trialgebra $(T, \ \prec, \ \succ, \ \circ)$. A {\it{$t$-Baxter operator over
the dendriform trialgebra}} $(T, \ \prec, \ \succ, \ \circ)$ is a linear map $\gamma: T \xrightarrow{} T$ such that for all $x,y \in T$,
\begin{eqnarray*}
\gamma(x) \succ \gamma(y) &=& \gamma(x \succ \gamma(y) + \gamma(x) \succ y + t x \succ y), \\
\gamma(x) \prec \gamma(y) &=& \gamma(x \prec \gamma(y) + \gamma(x) \prec y + t x \prec y), \\
\gamma(x) \circ \gamma(y)&=& \gamma(x \circ \gamma(y) + \gamma(x) \circ y + t x \circ y),
\end{eqnarray*}
so that, $\gamma(x) \star \gamma(y) = \gamma(x \star \gamma(y) + \gamma(x) \star y + t x \star y).$
\end{defi}
\Rk 
Observe that the dendriform trialgebra $(T, \ \prec, \ \succ, \ \circ)$ has another associative law $\circ$ and that $\gamma$ is a $t$-Baxter operator over $(T, \ \circ)$. Therefore, requiring the last axiom
means that $\gamma: (T, \ \star^{\circ} _{\gamma} ) \xrightarrow{} (T, \ \circ)$ is \footnote{We have placed $\circ$ in exponent to indicate that the associative product $\star^{\circ}$ comes from the associative product $\circ$.} a morphism of associative algebras and that $(T, \ \star^{\circ}_{\gamma})$ has also a dendriform trialgebra structure. 

Let $(T, \ \prec, \ \succ, \ \circ)$ be a dendriform trialgebra and $\gamma: T \xrightarrow{} T$ be a $t$-Baxter operator. For all $x,y \in T$, define new binary operations on $T$ by:
\begin{eqnarray*}
x \searrow_\gamma y &=& \gamma(x) \succ y, \ \
x \nearrow_\gamma y = x \succ \gamma(y), \ \
x \swarrow_\gamma y = \gamma(x) \prec y, \\
 x \nwarrow_\gamma y &=& x \prec \gamma(y), \ \
x \uparrow_\gamma y \ = x \circ \gamma(y), \ \ \ \ \
x \downarrow_\gamma y = \gamma(x) \circ y.
\end{eqnarray*}
Define also:
\begin{eqnarray*}
 x \succ_{\gamma} y &:=& x \searrow_{\gamma} y + x \nearrow_{\gamma} y + t x \succ y := \gamma(x) \succ y + x \succ \gamma(y) + t x \succ y, \\
 x \prec_{\gamma} y &:=& x \nwarrow_{\gamma} y + x \swarrow_{\gamma} y + t x \prec y := \gamma(x) \prec y + x \prec \gamma(y) + t x \prec y, \\
x \bar{\circ}_{\gamma} y &:=& x \uparrow_{\gamma} y + x \downarrow_{\gamma} y + t x \circ y :=  x \circ \gamma(y) + \gamma(x) \circ y + t x \circ y, \\
x \vee_{\gamma} y &:=& x \searrow_{\gamma} y + x \swarrow_{\gamma} y +  x \downarrow_{\gamma} y, \\
x \wedge_{\gamma} y &:=& x \nearrow_{\gamma} y + x \nwarrow_{\gamma} y +  x \uparrow_{\gamma} y.
\end{eqnarray*}
\begin{prop}
\label{hexd}
With these operations, the dendriform trialgebra $T$ has a $t$-ennea-algebra structure. As a consequence, 
the $k$-vector spaces $(T, \ \prec_\gamma, \ \succ_\gamma, \ \circ_\gamma)$ and $(T, \ \wedge_\gamma, \ \vee_\gamma, \ t \star)$ are dendriform trialgebras.
\end{prop}
\Proof
This is a consequence of Theorem \ref{enneathm}. Let $x,y \in T$. Observe that
$\gamma(x) \succ \gamma(y) = \gamma(x \succ_{\gamma} y) $,
$\gamma(x) \prec \gamma(y) = \gamma(x \prec_{\gamma} y) $ and 
$\gamma(x) \circ \gamma(y) = \gamma(x \circ_{\gamma} y) $, i.e., the map $\gamma$ is a morphism of dendriform trialgebras $(T,\ \prec_{\gamma},\ \succ_{\gamma},\ \circ_{\gamma}) \xrightarrow{} (T,\ \prec,\ \succ, \ \circ)$.
Moreover observe that
$x \vee_{\gamma} y := x \searrow_{\gamma} y + x \swarrow_{\gamma} y + x \downarrow_{\gamma} y:= \gamma(x) \succ y + \gamma(x) \prec y + \gamma(x) \circ y := \gamma(x) \star y$ and
$x \wedge_{\gamma} y := x \nearrow_{\gamma} y + x \nwarrow_{\gamma} y + x \uparrow_{\gamma} y:= x \succ \gamma(y) + x \prec \gamma(y) + x \circ \gamma(y) := x \star \gamma(y)$, which can be helpfull in computations. For instance, 
\begin{eqnarray*}
(x \nwarrow_{\gamma} y) \nwarrow_{\gamma} z &:=& (x \prec \gamma(y)) \prec \gamma (z), \\
&:=& x \prec (\gamma(y) \star \gamma (z)), \\
&:=& x \prec \gamma(y \star \gamma (z) + \gamma(y) \star z + tx \circ z), \\
&:=& x \nwarrow_{\gamma} (y \bar{\star}_{\gamma} z).
\end{eqnarray*}
\eproof
\begin{defi}{[Nested dendriform trialgebra]}
The $k$-vector space $(T, \ \prec,\ \succ,\ \circ)$ is said to be a {\it{nested dendriform trialgebra}} if it is a dendriform trialgebra and if the product of the associative algebra $(T,\ \circ)$ is a cluster of three operations verifying dendriform trialgebra axioms.
For instance, according to the previous Proposition,
under the action of a $t$-Baxter operator, $\gamma$, the horizontal and vertical dendriform trialgebras
$(T, \ \prec_\gamma,\  \succ_\gamma, \ \bar{\circ}_\gamma)$ and $(T, \ \wedge_\gamma, \ \vee_\gamma, \ t \star)$ have
their third (associative) laws, $\bar{\circ}_\gamma$ and $t \star$, which
both split into three operations turning $(T, \ \circ_\gamma)$ and $(T, \ t \star)$ into dendriform trialgebras.
\end{defi}
Let $(A, \ \cdot)$ be an associative algebra. To construct such $t$-ennea-algebras on $A$, we will consider two $t$-Baxter operators on $A$, $\beta$ and $\gamma$,  
which commute i.e., $\beta\gamma=\gamma\beta$. 
\begin{prop}
Let $\beta$ and $\gamma$ be a pair of commuting $t$-Baxter operators on an associative algebra $(A, \ \cdot)$. Then, $\gamma$ is a $t$-Baxter operator on the dendriform trialgebra $(A, \ \prec_{\beta}, \ \succ_{\beta},\ \circ_{\beta})$.
\end{prop}
\Proof
Let $x,y \in A$.
Let us check for instance the following equality:
$\gamma(x) \succ_{\beta} \gamma(y) = \beta(\gamma(x))\cdot \gamma(y)= \gamma(\beta(x))\cdot \gamma(y) := \gamma(\beta(x)\cdot \gamma(y) + \gamma(\beta(x))\cdot y + t \beta(x) \cdot y)=
\gamma(\beta(x)\cdot \gamma(y) + \beta(\gamma(x))\cdot y + t \beta(x) \cdot y)=
\gamma(x \succ_{\beta} \gamma(y) + \gamma(x)\succ_{\beta} y + t x\succ_{\beta} y).$
\eproof
\begin{coro}
\label{corohex}
Let $\beta$ and $\gamma$ be a pair of commuting $t$-Baxter operators on an associative algebra $(A, \ \cdot)$. Then, there is a $t$-ennea-algebra structure over the dendriform trialgebra $(A, \ \prec_{\beta}, \ \succ_{\beta},\ \circ_{\beta})$, with operations defined on $(A, \ \cdot)$ by:
$$x \searrow y = \beta(\gamma(x))\cdot y = \gamma(\beta(x))\cdot y, \  \
x \nearrow y = \beta(x)\cdot \gamma(y), \  \
x \swarrow y = \gamma(x)\cdot \beta(y),  \  \
$$
$$x \nwarrow y = x\cdot \beta(\gamma(y)) = x\cdot \gamma(\beta(y)), \ \
x \uparrow y = tx\cdot \gamma(y), \  \
x \downarrow y = t\gamma(x)\cdot y.$$
\end{coro}
\Proof
Apply Proposition \ref{hexd} to the $t$-Baxter operator
on the dendriform trialgebra $(A, \ \prec_{\beta}, \ \succ_{\beta}, \ \circ_{\beta})$.
\eproof  

\NB
As in \cite{AguiarLoday} (case $t=0$), if we start with applying first the 
$t$-Baxter operator $\gamma$ and then $\beta$, we will obtain another $t$-ennea-algebra over the dendriform trialgebra $(A, \ \prec_{\gamma}, \ \succ_{\gamma}, \ \circ_{\gamma})$, which is the transpose of the first one established in Corollary \ref{corohex}. 

\subsection{Relations with (left) pre-Lie algebras}
We now relate these structures to left pre-Lie algebras, called in the sequel pre-Lie algebras. 
A pre-Lie algebra is a $k$-vector space $P$ equipped with a binary operation $\bowtie: \ P \otimes P \xrightarrow{} P$ such that for all $x,y \in P$:
$$ x \bowtie (y \bowtie z) - (x \bowtie y ) \bowtie z = 
y \bowtie (x \bowtie z) - (y \bowtie x) \bowtie z.$$
In particular, any associative algebra is a (trivial) pre-Lie algebra. The new binary operation $P \otimes P \xrightarrow{} P$ defined by $[x,y]_{\bowtie}:= x \bowtie y - y \bowtie x $ turns $(P, [\cdot \ , \ \cdot]_{\bowtie})$ into a Lie algebra \cite{Gerst}.
\begin{prop}
\label{prelie}
Let $t \in k$. Let $(T, \ \prec, \ \succ, \ t\circ)$ be a dendriform trialgebra with $\star \longrightarrow \prec + \succ + t\circ$.
Consider the following operation: $x \bowtie y := x \succ y - y \prec x + tx \circ y$, constructed on the dendriform trialgebra $(T, \ \prec, \ \succ, \ t\circ)$. Then, $(T, \ \bowtie )$ is pre-Lie algebra. Moreover the two pre-Lie algebras $(T, \ \bowtie )$ and $(T, \ \star )$ give same Lie algebra, i.e., $[x,y]_{\bowtie} = [x,y]_{\star}$, for all $x,y \in T$.
\end{prop}
\Proof
Straightforward. Notice that for $t=0$, a dendriform dialgebra gives a pre-Lie structure on itself. 
\eproof
\begin{coro}
Let $t \in k$. Let $(T, \ \prec, \ \succ, \ \circ)$ be a dendriform trialgebra with $\star \longrightarrow \ \prec + \succ + \circ$ which admits a $t$-ennea-algebra structure, i.e., there exist on $T$, binary operations  $\triangleleft, \ \triangleright, \ \bar{\circ}, \ \wedge, \ \vee$, verifying Definition \ref{ennea} such that 
$\bar{\star} \longrightarrow \triangleright +\triangleleft + \bar{\circ}$ and $\bar{\star} \longrightarrow \wedge + \vee + t\star$. 
Then, there exist at least two pre-Lie structures on this dendriform trialgebra. That given by the following binary operation $x \bowtie y := x \triangleright y - y \triangleleft x + x \bar{\circ} y$ and that given by $x \ \hat{\bowtie} \ y := x \vee y - y \wedge x + tx \star y$. Moreover these two pre-Lie algebras $(T, \ \bowtie)$, $(T, \ \hat{\bowtie}) $ and $(T, \ \bar{\star})$ give the same Lie algebra, i.e., for all $x,y \in T$, $[x,y]_{\bowtie}=[x,y]_{\hat{\bowtie}}=[x,y]_{\bar{\star}}$.
\end{coro}
\Proof
Straightforward. 
\eproof
\NB
Observe that for $t=0$, a quadri-algebra gives two pre-Lie structures on itself which degenerate into the same Lie structure.
\NB
For a nested dendriform trialgebra, the third operation which is also a cluster of three operations
obeying dendriform trialgebra axioms gives rise to another pre-Lie structure according to Proposition \ref{prelie}.
As a $t$-ennea-algebra generates two nested dendriform trialgebras, we can say that at least four non-trivial pre-Lie algebras can be constructed on such a $k$-vector space.
\section{$t$-Ennea-algebras from $t$-infinitesimal bialgebras}
If $(A, \ \mu)$ is a $k$-algebra, then 
the following operations $a (x \otimes y) := ax \otimes y$ and $ (x \otimes y) a:= x \otimes ya$ defined for all $a, \ x, \ y \in A$
turn $A \otimes A$ into a $A$-bimodule.
Let us now yield a process, inspired from \cite{Aguiar, AguiarLoday} to produce commuting pair of $t$-Baxter operators.
\begin{defi}{[$t$-infinitesimal bialgebra]}
A {\it{$t$-infinitesimal bialgebra}} (abbreviated $\epsilon(t)$-bialgebra) is a triple $(A, \ \mu, \ \Delta )$ where $(A, \ \mu)$ is an associative $k$-algebra and $(A,\ \Delta )$ is a coassociative coalgebra such that
$ \Delta \mu := (\mu \otimes id)(id \otimes \Delta) + (id \otimes \mu)(\Delta \otimes id) + t \ id \otimes id $.
Therefore, for all $a,b \in A$,
$\Delta (ab)= a_{(1)} \otimes a_{(2)} b + a b_{(1)} \otimes b_{(2)}  + t a \otimes b$.
If $t=0$, a $t$-infinitesimal bialgebra is called
an {\it{infinitesimal bialgebra}} or a $\epsilon$-bialgebra. Such  bialgebras
appeared for the first time in the work of Joni and Rota in \cite{Rota}, see also Aguiar \cite{Aguiar}, for the case $t=0$ and Loday \cite{Lodayscd}, for the case $t=-1$.
\end{defi}

We now produce a pair of commuting $t$-Baxter operators.
The $k$-vector space $\textsf{End}(A)$ of linear
endomorphisms of $A$ is viewed as an associative algebra under composition denoted simply by concatenation $TS$, for $T,S \in \textsf{End}(A)$. In addition, equip $\textsf{End}(A)$ with the convolution product $*$, defined by $T *S := \mu(T \otimes S)\Delta$, for all $T,S \in \textsf{End}(A)$.
\begin{prop}
Let $(A, \ \mu, \ \Delta)$ be a $\epsilon(t)$-bialgebra and consider
$\textsf{End}(A)$ as an associative algebra under composition, equipped with the convolution product. There exist two commuting $t$-Baxter operators $\beta$ and $ \gamma$ on $\textsf{End}(A)$, defined for all $T \in \textsf{End}(A)$ by:
$$\beta(T):= id *T \ \textrm{and} \ \gamma(T):= T*id. $$
\end{prop}
\Proof
Let $a \in A$. Using Sweedler's notation, set 
$\Delta(a) := a_{(1)} \otimes a_{(2)}$. Let us check that
the map $\beta$ is a $t$-Baxter operator;
the verification for $\gamma$ is similar.
We have $\beta(S)(a)=a_{(1)}S(a_{(2)})$. Therefore,
$\Delta(\beta(S)(a))=a_{(1)}S(a_{(2)})_{(1)} \otimes S(a_{(2)})_{(2)} +
a_{(1)} \otimes a_{(2)(1)} S(a_{(2)(2)}) + t a_{(1)} \otimes S(a_{(2)}),$
thus
$\beta(T)\beta(S)(a) := a_{(1)}S(a_{(2)})_{(1)}T(S(a_{(2)})_{(2)}) +
a_{(1)} T(a_{(2)(1)} S(a_{(2)(2)})) + t a_{(1)} TS(a_{(2)})=
\beta(\beta(T)S + T\beta(S) + tTS)(a).$
\eproof
\begin{coro}
\label{zeze}
Let $(A,\ \mu, \ \Delta )$ be a $\epsilon(t)$-bialgebra. There is a $t$-ennea-algebra structure over the dendriform trialgebra (\textsf{End}$(A), \ \prec_{\beta}, \ \succ_{\beta}, \ \circ_{\beta}$) defined by:
$$T  \searrow S=(id*T*id)S; \ \ \
T  \nearrow S=(id*T)(S*id); \ \ \
T  \swarrow S=(T*id)(id*S); \ \ \
T  \nwarrow S=T(id*S*id). $$
$$T  \downarrow S= t(T*id)S; \ \ \
T  \uparrow S= tT(S*id).$$
\end{coro}
\Proof
Apply the Corollary \ref{corohex} to the pair of commuting $t$-Baxter operators $\beta$ and $\gamma$.
\eproof

\noindent
Recall there exists also a $t$-ennea-algebra structure over the dendriform trialgebra (\textsf{End}$(A), \ \prec_{\gamma}, \ \succ_{\gamma}, \ \circ_{\gamma}$) which is the transpose of that established in Corollary \ref{zeze}.
\section{Pre-Lie algebras and $t$-infinitesimal bialgebras}
To generalise \cite{AguiarLoday}, the notion of $t$-infinitesimal bialgebras, $t \in k$, has been introduced. Let us take the opportunity of this paper to show relations between $t$-infinitesimal bialgebras and pre-Lie algebras, relations already established in \cite{Aguiar} for $t=0$.
For that, we show that any $\epsilon(t)$-bialgebra inherits a (left) pre-Lie structure. Following M. Aguiar \cite{Aguiar}, we construct
a (left) pre-Lie structure on a $\epsilon(t)$-bialgebra.
\begin{prop}
Let $t \in k$.
Let $(A, \ \mu, \ \Delta)$ be a $\epsilon(t)$-bialgebra and $a,b \in A$. Using Sweedler notation, we set $\Delta(b):=b_{(1)} \otimes b_{(2)}$. Then, the new operation on $A$ defined by:
$$ a \bowtie b := b_{(1)} a b_{(2)},$$
turns $A$ into a (left)-pre-Lie algebra, i.e., the law $\bowtie$ verifies for all $x,y \in A$:
$$ x \bowtie (y \bowtie z) - (x\bowtie y) \bowtie z=
y \bowtie (x \bowtie z) - (y \bowtie x) \bowtie z.$$ 
Therefore $(A, \ [\cdot, \ \cdot]_{\bowtie})$, where $[x,y]_{\bowtie} := x \bowtie y -y \bowtie x$, for all $x,y \in A$, is a Lie algebra.
\end{prop}
\Proof
Straightforward.
\eproof

\noindent
Let $(A, \ \mu, \ \Delta)$ be a $\epsilon(t)$-bialgebra. 
The space $\textsf{End}(A)$ can be viewed as a Lie algebra
under the commutator bracket $[T, S]:= TS -ST$.
Denote by $\textsf{Der}(A, \mu)$, the $k$-vector space of 
derivations $D: A \xrightarrow{} A$ on the associative algebra $(A, \mu)$. 
\begin{prop}
Let $(A, \ \mu, \ \Delta)$ be a $\epsilon(t)$-bialgebra, where $t \in k$. 
Then, the operator $L: a \mapsto L_a: x \mapsto a \bowtie x + tax$ is defined on $A$ with values in $\textsf{Der}(A, \mu)$. 
\end{prop}
\Proof
Let $(A, \ \mu, \ \Delta)$ be a $\epsilon(t)$-bialgebra, where $t \in k$. 
Fix $a \in A$. It is easy to show that $L_a(x  y ) = xL_a( y ) + L_a(x )y$, for all $x,y \in A$. 
\eproof
\begin{defi}{}
Let $(A, \ \mu, \ \Delta)$ be a $\epsilon(t)$-bialgebra. 
Similarly, following \cite{Aguiar}, we define
a {\it{bi-derivation}} on $A$ as a map $B: A \xrightarrow{} A$ which is both a derivation on $(A, \ \mu)$, i.e.,
$B(ab):= B(a)b + aB(b)$ and a coderivation on $(A, \ \Delta)$, i.e., $\Delta(B(b)):= b_{(1)} \otimes B(b_{(2)}) +
B(b_{(1)}) \otimes b_{(2)}$, for all $a,b \in A$.
\end{defi}
\begin{prop}
Let $(A, \ \mu, \ \Delta)$ be a $\epsilon(t)$-bialgebra and  $B: A \xrightarrow{} A$ be a bi-derivation on $A$. Then, $B$ is a derivation on the pre-Lie algebra $(A, \ \bowtie)$.
\end{prop}
\Proof
It suffices to show that $B(a \bowtie b)= B(a) \bowtie b + B(a) \bowtie b$, for all $a,b \in A$, which is straightforward.
\eproof
\section{Examples}
\label{examp}
In Section~\ref{examp}, we point out a link between (weighted) graph theory and our work.

\subsection{Path algebras of weighted directed graphs and quadri-algebras}
\label{grappph}
Let $G$ be a weighted directed graph and $G_n$ be the set of path of length $n$, in particular $G_0$ is the vertex set  and $G_1$ the arrow set. A directed arrow $e_i \longrightarrow e_j$, with $e_i, \ e_j \in G_0$, will be denoted by $a_{(i,j)}$. A directed graph is said to be {\it{weighted}} if it is equipped with a weight map, i.e., with a map $w$ from $G_1$ to $k$.
The relations
$e_ie_j =1$ if $i=j$ and 0 otherwise, $a_{(i,j)} a_{(k,l)} = 0$ if $j \not= k$, $e_ia_{(k,l)} = a_{(k,l)}$ if $i=k$ and 0 otherwise and $a_{(k,l)}e_i = a_{(k,l)}$ if $i=l$ and 0 otherwise, for every $e_i, \ e_j \in G_0$ and $a_{(i,j)}, \ a_{(k,l)} \in G_1$, turn the $k$-vector space $kG := \oplus_{n=0} ^{ \infty} kG_n$ into an associative algebra called the path algebra. The product is then the concatenation of paths whenever possible. Let us define the following co-operation $\Delta$, for any $e_i \in G_0$ by $\Delta(e_i):=0$, for any weighted arrow $a_{(i,j)}$ by $\Delta a_{(i,j)}:= w(a_{(i,j)}) e_i \otimes e_j$ and for any weighted path $\alpha:= a_{(i_1,i_2)} a_{(i_2,i_3)} \ldots a_{(i_{n-1},i_n)}$ by, 
$$\Delta (\alpha):= \Delta (a_{(i_1,i_2)}) a_{(i_2,i_3)} \ldots a_{(i_{n-1},i_n)} + \ldots + a_{(i_1,i_2)} a_{(i_2,i_3)} \ldots a_{(i_{n-2},i_{n-1})}\Delta a_{(i_{n-1},i_n)}.$$

\begin{prop}
\label{al}
Equipped with this co-operation, any weighted directed graph carries a $\epsilon$-bialgebra structure.
\end{prop}
\Proof
Straightforward.
\eproof
\NB
We recover the $\epsilon$-bialgebra structure found by Aguiar \cite{Aguiar} in the case of the trivial weight map, i.e., $w(a_{(i,j)})=1$ for all $a_{(i,j)} \in G_1$. It can be also interesting to observe that weight maps are used to study random walks on directed graphs. The conditions on the weight map are the following: $w: G_1 \xrightarrow{} \mathbb{R}_+$ and for any vertex $e_i \in G_0$, $\sum_{j} w(a_{(i,j)}) =1$.

As a corollary, the path algebra of a weighted directed graph can be equipped with a pre-Lie algebra
and the $k$-vector space $\textsf{End}(kG)$ has a structure of quadri-algebra ($t$-ennea-algebra with $t=0$). 
\begin{prop}
Keep notation introduced in this subsection.
Let $G$ be a directed graph and $kG$ its path algebra. Let $w^1$ and $w^2$ be two weight maps. Define two coproducts $\Delta_1, \ \Delta_2: kG \xrightarrow{} kG^{\otimes 2}$ by 
$\Delta_2(e_i) := \Delta_1(e_i):=0$ for all $e_i \in G_0$, for all arrows $a_{(i,j)}$ by $\Delta_1 a_{(i,j)}:= w^1(a_{(i,j)}) e_i \otimes e_j$ and 
$\Delta_2 a_{(i,j)}:= w^2(a_{(i,j)}) e_i \otimes e_j$.
For any path $\alpha:= a_{(i_1,i_2)} a_{(i_2,i_3)} \ldots a_{(i_{n-1},i_n)}$ by, 
$$\Delta_p (\alpha):= \Delta_p (a_{(i_1,i_2)}) a_{(i_2,i_3)} \ldots a_{(i_{n-1},i_n)} + \ldots + a_{(i_1,i_2)} a_{(i_2,i_3)} \ldots a_{(i_{n-2},i_{n-1})}\Delta_p a_{(i_{n-1},i_n)},$$
for $p=1,2$.
Then, for all $p,q=1,2$,
$$ (id \otimes \Delta_p)\Delta_q = (\Delta_q \otimes id)\Delta_p. $$
\end{prop}
\Proof
Straightforward.
\eproof
\subsection{Path algebras of directed graphs and $(-1)$-ennea-algebras}
Notation of Subsection \ref{grappph} are kept.
We consider a directed graph and its path algebra (the graph is no longer weighted or equipped with the trivial weight map). 
Let $\alpha:= a_{(i_1,i_2)} a_{(i_2,i_3)} \ldots a_{(i_{n-1},i_n)}$ be a path. Define the co-operation $\hat{\Delta}$ for any path as follows: $\hat{\Delta} \alpha:= s(\alpha) \otimes \alpha + a_{(i_1,i_2)} \otimes a_{(i_2,i_3)} \ldots a_{(i_{n-1},i_n)} + \ldots + \alpha \otimes  t(\alpha)$, where $s$ is the source map and $t$, the terminus map and for any vertices $e_i \in G_0$, $\hat{\Delta} e_i:= e_i \otimes e_i$.
\begin{prop}
\label{dup}
Let $G$ be a directed graph.
Equipped with the co-operation $\hat{\Delta}$, the path algebra of $G$ is a $\epsilon(-1)$-bialgebra.
\end{prop}
\Proof
Straightforward.
\eproof
\begin{coro}
Fix an integer $n$, possibly equal to infinity. Let $S=\{X_1, \ldots, \ X_n \}$ be a set. Let $k \bra S \ket$ be the non-commutative polynomials constructed from $X_1, \ldots, \ X_n$. Define the co-operation
$\Delta: k \bra S \ket \xrightarrow{} k \bra S \ket^{\otimes 2} $ such that $\Delta(X_i):= X_i \otimes 1 + 1 \otimes X_i$, $\Delta (1) := 1 \otimes 1$ and
$\Delta(X_{i_1}X_{i_2}X_{i_3} \ldots X_{i_p})
= 1 \otimes X_{i_1} X_{i_2}X_ {i_3} \ldots X_{i_p} +
X_{i_1} \otimes X_{i_2}X_ {i_3} \ldots X_{i_p}
+ X_{i_1} X_{i_2} \otimes X_ {i_3} \ldots X_{i_p}
+ \ldots + X_{i_1} X_{i_2}X_ {i_3} \ldots X_{i_p} \otimes 1,$ for all $X_{i_1}\ldots X_{i_p} \in S$.
Then, $\Delta$ is a coassociative coproduct and the $k$-vector space $(k \bra S \ket, \ \Delta)$ is  a $\epsilon(-1)$-bialgebra.
\end{coro}
\Proof
It suffices to consider the graph with one vertex and $n$ loops.
\eproof

As another corollary, the path algebra of a directed graph can be also equipped with another pre-Lie algebra. In addition, the $k$-vector space $\textsf{End}(kG)$ is a $(-1)$-ennea-algebra.
\subsection{Other constructions} 
A $t$-ennea-algebra can be constructed from any tensor product of two dendriform trialgebras.
\begin{prop}
\label{tensor}
Fix $t \in k$.
Let $(A,\ \prec, \ \succ, \ \circ)$ and $(B,\ \prec', \ \succ', \ t\circ') $ be two dendriform trialgebras. Set $\star \longrightarrow \ \prec + \succ + \circ$ and  $\star' \longrightarrow \ \prec' + \succ' + t\circ'$.
Define for all $a_1,a_2 \in A$ and $b_1, b_2 \in B$ the following operations:
\begin{eqnarray*}
& &(a_1 \otimes a_2) \nwarrow (b_1 \otimes b_2) := (a_1 \prec a_2) \otimes (b_1 \prec' b_2); \ \ \ \ \ \
(a_1 \otimes a_2) \swarrow (b_1 \otimes b_2) := (a_1 \prec a_2) \otimes (b_1 \succ' b_2); \\
& &(a_1 \otimes a_2) \nearrow (b_1 \otimes b_2) := (a_1 \succ a_2) \otimes (b_1 \prec' b_2); \ \ \ \ \ \
(a_1 \otimes a_2) \searrow (b_1 \otimes b_2) := (a_1 \succ a_2) \otimes (b_1 \succ' b_2); \\ 
& &(a_1 \otimes a_2) \uparrow (b_1 \otimes b_2) := (a_1 \circ a_2) \otimes (b_1 \prec' b_2); \ \ \ \ \ \ \ \ 
(a_1 \otimes a_2) \downarrow (b_1 \otimes b_2) := (a_1 \circ a_2) \otimes (b_1 \succ' b_2);\\ 
& &(a_1 \otimes a_2) [\circ] (b_1 \otimes b_2) := (a_1 \circ a_2) \otimes (b_1 \circ' b_2); \ \ \ \ \ \ \ \ \
(a_1 \otimes a_2) [\prec] (b_1 \otimes b_2) := (a_1 \prec a_2) \otimes (b_1 \circ' b_2); \\
& &(a_1 \otimes a_2) [\succ] (b_1 \otimes b_2) := (a_1 \succ a_2) \otimes (b_1 \circ' b_2). \ \ \ \
\end{eqnarray*}
Then, the $k$-vector space $A \otimes B$ equipped with these operations is a $t$-ennea-algebra.
The vertical structure is given by:
\begin{eqnarray*}
(a_1 \otimes a_2) \wedge (b_1 \otimes b_2) &:=& (a_1 \star a_2) \otimes (b_1 \prec' b_2), \\
(a_1 \otimes a_2) \vee (b_1 \otimes b_2) &:=& (a_1 \star a_2) \otimes (b_1 \succ' b_2), \\
(a_1 \otimes a_2) [t\star] (b_1 \otimes b_2) &:=& t(a_1 \star a_2) \otimes (b_1 \circ' b_2).
\end{eqnarray*}
The horizontal structure is given by:
\begin{eqnarray*}
(a_1 \otimes a_2) \triangleright (b_1 \otimes b_2) &:=& (a_1 \succ a_2) \otimes (b_1 \star' b_2), \\
(a_1 \otimes a_2) \triangleleft (b_1 \otimes b_2) &:=& (a_1 \prec a_2) \otimes (b_1 \star' b_2), \\
(a_1 \otimes a_2) \bar{\circ} (b_1 \otimes b_2) &:=& (a_1 \circ a_2) \otimes (b_1 \star' b_2).
\end{eqnarray*}
The associative structure is then the expected one since, 
$$(a_1 \otimes a_2) \bar{\star} (b_1 \otimes b_2) \ :=  \ (a_1 \star a_2) \otimes (b_1 \star' b_2). $$
\end{prop}
\Proof
Straightforward.
\eproof
\begin{prop}
Let $S$ be a set and $\Delta: kS \xrightarrow{} kS^{\otimes 2}$ a coassociative coproduct. Fix $t_1, t_2 \in k$.
Consider the free $k$-algebra $As(S)$ generated by $S$. Extend $\Delta$ by $\Delta^{t_i}$, $i=1,2$ as follows. On $kS$, $\Delta^{t_i}= \Delta$ and for all $u,v \in As(S)$, 
$$\Delta^{t_i} (uv) := \Delta^{t_i}(u)v + u \Delta^{t_i}(v) + t_i u \otimes v.$$
Then for $i=1,2$, $(As(S), \ \Delta^{t_i})$ is a $\epsilon(t_i)$-bialgebra. Moreover, for all $i,j=1,2$,
$$ (id \otimes \Delta^{t_i})\Delta^{t_j} = (\Delta^{t_j} \otimes id)\Delta^{t_i}. $$
\end{prop}
\Proof
Straightforward.
\eproof
\section{Free $t$-ennea-algebra and open problems}
\label{freee}
In Section \ref{freee}, we discuss some properties of the free $t$-ennea-algebra, relations with combinatorics and connected Hopf algebras. Before, let us recall what an operad is, see \cite{Lodayscd, Fresse, GK} for instance. 

\noindent
Let $P$ be a type of algebras, for instance the $t$-ennea-algebras, and $P(V)$ be the free $P$-algebra on the $k$-vector space $V$. Suppose $P(V) := \oplus_{n \geq 1} P(n) \otimes _{S_n} V^{\otimes n},$ where $P(n)$ are right $S_n$-modules. Consider $P$ as an endofunctor on the category of $k$-vector spaces. The structure of the free $P$-algebra of $P(V)$ induces a natural transformation $\pi: P \tilde{\circ} P \xrightarrow{} P$ as well as $u: Id \xrightarrow{} P$ verifying usual associativity and unitarity axioms. An algebraic operad is then a triple $(P, \ \pi, \ u)$. A $P$-algebra is then a $k$-vector space $V$ together with a linear map $\pi_A: P(A) \xrightarrow{} A$ such that $\pi_A \tilde{\circ} \pi(A) = \pi_A \tilde{\circ} P(\pi_A)$ and $\pi_A \tilde{\circ} u(A) = Id_A.$
The $k$-vector space $P(n)$ is the space of $n$-ary operations for $P$-algebras. We will always suppose there is, up to homotheties, a unique $1$-ary operation, the identity, i.e.,  $P(1):= kId$ and that all possible operations are generated by composition from $P(2)$. The operad is said to be {\it{binary}}. It is said to be {\it{quadratic}} if all the relations between operations are consequences of relations described exclusively with the help of monomials with two operations. An operad is said to be {\it{non-symmetric}} if, in the relations, the variables $x,y,z$ appear in the same order. The $k$-vector space $P(n)$ can be written as  $P(n):= P'(n) \otimes k[S_n]$, where $P'(n)$ is also a $k$-vector space and $S_n$ the symmetric group on $n$ elements.
In this case, the free $P$-algebra is entirely induced by the free $P$-algebra on one generator $P(k):= \oplus_{n \geq 1} \ P'(n)$. The generating function of the operad $P$ is given by:
$$ f^{P}(x):= \sum \ (-1)^n \frac{\textrm{dim} \ P(n)}{n!} x^n := \sum \ (-1)^n \textrm{dim} \ P'(n) x^n.$$
Below, we will indicate the sequence $(\textrm{dim} \ P'(n))_{n \geq 1}$.

Let $V$ be a $k$-vector space. The {\it{free $t$-ennea-algebra}} $\mathcal{E}(V)$ on $V$ is by definition,
a $t$-ennea-algebra equipped with a map $i: \ V \mapsto \mathcal{E}(V)$ which satisfies the following universal property:
for any linear map $f: V \xrightarrow{} A$, where $A$ is a $t$-ennea-algebra, there exists a unique $t$-ennea-algebra morphism $\bar{f}: \mathcal{E}(V) \xrightarrow{} A$ such that $\bar{f} \circ i=f$.

Since the nine operations of a $t$-ennea-algebra have no symmetry and since compatibility axioms involve only monomials where $x, \ y$ and $z$ stay in the same order, the free $t$-ennea-algebra is of the form:
$$ \mathcal{E}(V) := \bigoplus_{n \geq 1} \ \mathcal{E}_n \otimes V^{\otimes n} .$$
In particular, the free $t$-ennea-algebra on one generator $x$ is $ \mathcal{E}(k) := \bigoplus_{n \geq 1} \mathcal{E}_n,$ where $\mathcal{E}_1:= kx$, $\mathcal{E}_2:= k(x\uparrow x) \oplus k(x\downarrow x)\oplus k(x\searrow x) \oplus k(x \nearrow x)\oplus k(x\swarrow x) \oplus k( x\nwarrow x) \oplus 
k(x \prec x) \oplus  k(x \succ x) \oplus k( x \circ x).$
The space of three variables made out of nine operations is of dimension $2 \times 9^2= 162$. As we have $7 \times 7$ relations, the space $\mathcal{E}_3$ has a dimension equal to $162-49=113.$
Therefore, the sequence associated with the dimensions of $(\mathcal{E}_n)_{n \in \mathbb{N}}$ starts with $1, \ 9, \ 113 \ldots$
Finding the free $t$-ennea-algebra on one generator is an open problem.  For $t =0$, this object is called the free quadri-algebra \cite{AguiarLoday}. The free quadri-algebra $\mathcal{Q}(V)$ on a $k$-vector space $V$ is of the form, $$ \mathcal{Q}(V) := \bigoplus_{n \geq 1} \ \mathcal{Q}_n \otimes V^{\otimes n} .$$
On one generator, one finds the dimension of $\mathcal{Q}_1$ is equal to 1, of $\mathcal{Q}_2$ is equal to 4 and of $\mathcal{Q}_3$ is equal to 23.

Quadri-algebras give dendriform dialgebras. It is known that
the free dendriform dialgebra on one generator is related to planar binary trees \cite{Loday}. It is conjectured in \cite{AguiarLoday} that the degree $n$ part of the free quadri-algebra on one generator can be indexed by the non-crossing connected graphs with $(n+1)$ vertices. A hint is to observe that the sequence $1,4,23, \ldots$ starts like the sequence $d_1=1,\ d_2=4,\ d_3=23,\ d_4=156,\ d_5=1162,\ldots$, where $d_n$ counts the number of non-crossing connected graphs made with $n+1$ vertices \cite{Domb,Flajolet}. 
Similarly, ennea-algebras give dendriform trialgebras. It is known that the free dendriform trialgebra on one generator is related to planar rooted trees \cite{LodayRonco}. We do not know which object could play the r\^ole of non-crossing connected graphs in the case of the free ennea-algebra on one generator since unfortunately, for the time being, no sequence starting with $1,9,113, \ldots$ seems to be related to known combinatorics. 

However, a little bit more on the free 
$t$-ennea-algebra, for all $t \in k$, can be said. In fact, we will show that the augmented free $t$-ennea-algebra over a $k$-vector space $V$ has a connected Hopf algebra structure. Before, some preparations are needed. To be as self-contained as possible, we introduce some notation to expose a theorem due to Loday \cite{Lodayscd}. 

\noindent
Recall that a bialgebra $(H, \mu, \ \Delta, \ \eta, \ \kappa)$
is a unital associative algebra $(H, \mu, \ \eta)$ together with co-unital coassociative coalgebra $(H, \ \Delta, \ \kappa)$. Moreover, it is required that the coproduct $\Delta$
and the counit $\kappa$ are morphisms of unital algebras. A bialgebra is connected if there exists a filtration $(F_rH)_r$ such that $H = \bigcup_r F_rH$, where $F_0H := k1_H$ and for all $r$,
$$ F_rH := \{x \in H; \ \Delta(x) -1_H \otimes x - x\otimes 1_H \in F_{r-1}H \otimes F_{r-1}H \}.$$
Such a bialgebra admits an antipode. Consequently, connected bialgebras are connected Hopf algebras.  

\noindent
Let $P$ be a binary quadratic operad. By a \textit{unit action} \cite{Lodayscd}, we mean the choice of two linear applications:
$$\upsilon: P(2) \xrightarrow{} P(1), \ \ \ \ \ \varpi:P(2) \xrightarrow{} P(1),$$
giving sense, when possible, to $x \circ 1$ and $1 \circ x$, for all operations $\circ \in P(2)$ and for all $x$ in the $P$-algebra $A$, i.e.,
$x \circ 1 = \upsilon(\circ)(x)$ and $1 \circ x= \varpi(\circ)(x)$.
If $P(2)$ contains an associative operation, say $\bar{\star}$, then we require that $x \bar{\star} 1 := x := 1 \bar{\star} x$, i.e., $\upsilon(\bar{\star}) := Id := \varpi(\bar{\star})$.
We say that the unit action, or the couple $(\upsilon,\varpi)$ is {\it{compatible}} with the relations of the $P$-algebra $A$ if they still hold on $A_+:= k1 \oplus A$ as far as the terms as defined.
Let $A$, $B$ be two $P$-algebras such that $P(2)$ contains an associative operation $\bar{\star}$. Using the couple $(\upsilon,\varpi)$, we extend binary operations $\circ \in P(2)$ to the $k$-vector space $A \otimes 1.k \oplus k.1 \otimes B \oplus A \otimes B$ by requiring:
\begin{eqnarray}
(a \otimes b) \circ (a' \otimes b') &:= &(a \bar{\star} a') \otimes (b \circ b') \ \ \ \textrm{if} \ \ b \otimes b' \not= 1 \otimes 1, \\
(a \otimes 1) \circ (a' \otimes 1) &:= &(a \circ a') \otimes 1, \ \ \ \textrm{otherwise}.
\end{eqnarray}
The unit action or the couple $(\upsilon,\varpi)$ is said to be {\it{coherent}} with the relations of $P$ if $A \otimes 1.k \oplus k.1 \otimes B \oplus A \otimes B$, equipped with these operations is still a $P$-algebra. Observe that a necessary condition for having coherence is compatibility.

\noindent
One of the main interest of these two concepts is the construction of a connected Hopf algebra on the augmented free $P$-algebra.
\begin{theo} [\textbf{Loday \cite{Lodayscd}}]
\label{Lodaych}
Let $P$ be a binary quadratic operad. Suppose there exists an associative operation in $P(2)$. Then, any unit action coherent with the relations of $P$ equips the augmented free $P$-algebra $P(V)_+$ on a $k$-vector space $V$ with a
coassociative coproduct $\Delta: P(V)_+ \xrightarrow{}  P(V)_+ \otimes P(V)_+,$
which is a $P$-algebra morphism.
Moreover,
$P(V)_+$ is a connected Hopf algebra.
\end{theo}
\Proof
See \cite{Lodayscd} for the proof. However, we reproduce it to make things clearer.
Let $V$ be a $k$-vector space and $P(V)$ be the free $P$-algebra on $V$.
Since the unit action is coherent, $P(V)_+ \otimes P(V)_+$ is a $P$-algebra. Consider the linear map $\delta: V \xrightarrow{} P(V)_+ \otimes P(V)_+$, given by $v \mapsto 1 \otimes v + v \otimes 1$. Since $P(V)$ is the free $P$-algebra on $V$, there exists a unique extension of $\delta$ to a morphism of augmented $P$-algebra $\Delta: P(V)_+ \xrightarrow{}  P(V)_+ \otimes P(V)_+$. Now, $\Delta$ is coassociative since the morphisms $(\Delta \otimes id)\Delta$ and 
$(id \otimes \Delta)\Delta$ extend the linear map $ V \xrightarrow{} P(V)_+ ^{\otimes 3}$ which maps $v$ to $1 \otimes 1 \otimes v + 1 \otimes v \otimes 1 + v \otimes 1 \otimes 1$. By unicity of the extension, the coproduct $\Delta$ is coassociative. The bialgebra so obtained is connected. Indeed, by definition, the free $P$-algebra $P(V)$ can be written as
$P(V) := \oplus_{n \geq 1} \ P(V)_n$, where $P(V)_n$ is the $k$-vector space of products of $n$ elements of $V$. Moreover, we have $\Delta(x \circ y ) := 1 \otimes (x \circ y )  + (x \circ y) \otimes 1 + x \otimes (1 \circ y) + y \otimes (x \circ 1)$, for all $x,y \in P(V)$ and $\circ \in P(2)$. The filtration of $P(V)_+$ is then $FrP(V)_+ = k.1 \oplus \bigoplus_{1 \leq n \leq r} P(V)_n$.
Therefore, $P(V)_+:=\cup_r \ FrP(V)_+$ and $P(V)_+$ is a connected bialgebra.
\eproof

\noindent
We will use this theorem to show that there exist a connected Hopf algebra structure on the augmented free $t$-ennea-algebra as well as on the augmented free commutative $t$-ennea-algebra, for all $t \in k$.
\begin{prop}
Let $\mathcal{E}(V)$ be the free $t$-ennea-algebra on a $k$-vector space $V$. Extend the binary operations $\nwarrow$ and $\searrow$ to $\mathcal{E}(V)_+$ as follows:
$$ e \nwarrow 1 := e, \  \ 1 \nwarrow e := 0, \ \ 1 \searrow e := e, \ \ e \searrow 1 := 0, \ \ \forall e \in \mathcal{E}(V).$$
In addition, for all operations $\diamond \in \mathcal{E}(2)$ different from $\nwarrow$ and $\searrow$ choose: 
$ e \diamond 1 := 0, \  \ 1 \diamond e := 0, \ \ \forall e \in \mathcal{E}(V).$ Then, for all $e \in \mathcal{E}(V)$:
$$ e \triangleleft 1 =e, \ \ 1 \triangleright e =e, \ \ 1 \vee e := e , \ \ e \wedge 1 := e, \ \ e \bar{\star} 1 = e = 1 \bar{\star} e.$$ 
Moreover, this choice is coherent.
\end{prop}
\Rk
We cannot extend the operations $\searrow$ and $\nwarrow$ to $k$, i.e., $1 \searrow 1$ and $ 1 \nwarrow 1$ are not defined.

\Proof
Keep notation introduced in that section.
Firstly, let us show that this choice is compatible. Let $x,y,z \in \mathcal{E}(V)_+$. We have to show for instance that the relation $(x \nwarrow y)\nwarrow z = x \nwarrow (y \bar{\star} z)$ holds in $\mathcal{E}(V)_+$. Indeed, for $x=1$ we get $0=0$. For $y =1$ we get $ x \nwarrow z = x \nwarrow z$ and for $z=1$ we get $ x \nwarrow y = x \nwarrow y$. We do the same thing with the 48 others and quickly found that the augmented free $t$-ennea-algebra $\mathcal{E}(V)_+$ is still a $t$-ennea-algebra, $t \in k$.

\noindent
Secondly, let us show that this choice is coherent. Let $x_1,x_2, x_3, y_1, y_2, y_3 \in \mathcal{E}(V)_+$. We have to show that, for instance:
$$(Eq. \ 11) \ \ \  ((x_1 \otimes y_1) \nwarrow (x_2 \otimes y_2)) \nwarrow (x_3 \otimes y_3) = (x_1 \otimes y_1) \nwarrow ((x_2 \otimes y_2) \bar{\star} (x_3 \otimes y_3)).$$
Indeed, if there exists a unique $y_i =1$, the other belonging to $\mathcal{E}(V)$, then, by definition we get:
$$ (x_1 \bar{\star} x_2 \bar{\star} x_3 ) \otimes (y_1 \nwarrow y_2) \nwarrow  y_3 = (x_1 \bar{\star} x_2 \bar{\star} x_3 ) \otimes y_1 \nwarrow ( y_2 \bar{\star}  y_3),$$
which always holds since our choice of the unit action is compatible. Similarly if $y_1=y_2=y_3=1$.
If $y_1=y_2=1$ and $y_3 \in \mathcal{E}(V)$, we get: $0=0$, similarly if $y_1=1=y_3$ and $y_2 \in \mathcal{E}(V)$. If $y_1 \in \mathcal{E}(V)$ and $y_2=1=y_3$, the two hand sides of $(Eq. \ 11)$ are equal to $(x_1 \bar{\star} x_2 \bar{\star} x_3 ) \otimes y_1$. Therefore, $ (Eq. \ 11)$ holds in $\mathcal{E}(V)  \otimes 1.k \oplus k.1 \otimes \mathcal{E}(V)  \oplus \mathcal{E}(V)  \otimes \mathcal{E}(V) $. Checking the same thing with the 48 other relations shows that our choice of the unit action is coherent.
\eproof
\begin{coro}
There exists a connected Hopf algebra structure on the augmented free $t$-ennea-algebra as well as on the augmented free commutative $t$-ennea-algebra, for all $t \in k$.
\end{coro}
\Proof
The first claim comes from the fact that our choice is coherent and from Theorem \ref{Lodaych}. For the second remark, observe that our choice is in agreement with the symmetry relations defining a commutative $t$-ennea-algebra since for instance $e \nwarrow^{op} 1 := 1 \searrow e :=e$ and $1 \searrow^{op} e := e \nwarrow 1 :=e,$ for all $e \in \mathcal{E}(V)$.
\eproof
\Rk
Such structures were found by Loday for quadri-algebras \cite{Lodayscd}.
\section{Baxter operators and formal deformations}
We end this article by showing the power of Baxter operators for the study of formal deformations of dendriform di- and tri-algebras, quadri-algebras and ennea-algebras. 

\subsection{1-parameter formal deformations} 
Recall \cite{Gerst1} that a formal deformation of a $k$-algebra $(A, \ \bullet)$ is a 1-parameter family of binary operations $(\bullet_i)_{i \geq 0}$ obtained by perturbing that of $A$, i.e.,
$$\bullet(h) \longrightarrow \ \bullet_0 + \bullet_1h + \bullet_2h^2 + \bullet_3h^3 + \ldots,$$ 
with $\bullet_0 := \bullet$. 
For this to make sense, we have to extend the coefficients to the ring $k[[h]]$ and consider $\bullet(h)$ as the binary operation $\bullet(h): A[[h]] \otimes_{k[[h]]} A[[h]] \xrightarrow{} A[[h]].$
In the category of associative algebras, it is required that $\bullet(h)$ has to be associative which entails conditions on the operations $\bullet_i$. For instance, $\bullet_1$ has to be a 2-Hochschild cocycle \cite{Gerst1}, i.e., has to verify,
$$ (x \bullet_1 y)\bullet z + (x \bullet y)\bullet_1 z := x \bullet_1 (y\bullet z) + x \bullet (y\bullet_1 z), \ \ \forall x,y,z \in A.$$
In the sequel, by a trivial formal deformation of an algebra of type $P$, say $A$, we mean the algebra of type $P$, $A[[h]]$, such that $x \diamond(h) y := x \tilde{\diamond} y,$ for all $\diamond \in P(2)$ where
$\tilde{\diamond}: A[[h]] \otimes_{k[[h]]} A[[h]] \xrightarrow{} A[[h]]$ extends the operation $\diamond$.

In the sequel, we will be interested in deforming special associative algebras such as dendriform di- and tri-algebras, quadri-algebras and ennea-algebras.
Let $\mathcal{A}$ and $\mathcal{B}$ be two small categories such that $\mathcal{A}$ is the category of $k$-algebras of type $P_\mathcal{A}$ and 
$\mathcal{B}$ the category of $k$-algebras of type $P_\mathcal{B}$. Suppose
$\mathcal{A} \subset \mathcal{B}$.
Therefore,
algebras of type $P_\mathcal{A}$ are also algebras of type $P_\mathcal{B}$. 
A
{\it{1-parameter formal deformation}} of a $k$-algebra $A$ of type $P_\mathcal{A}$ is the data of a formal deformation of $A$, i.e. $(A[[h]], \{\diamond_i(h) \}_{i \in I})$ and a family of $k$-algebras of type $P_\mathcal{B}$ denoted by $(B_\tau[[h]]:=(A[[h]], \{\diamond^\tau _i(h) \}_{i \in J})_{\tau \in k}$ with $I \subseteq  J$ and $B_0[[h]]:=(A[[h]], \{\diamond^0 _i(h) \}_{i \in J})$ equals to $(A[[h]], \{\diamond_i(h) \}_{i \in I})$. 

Let us show what this means in the category $\textsf{Ennea}$ of ennea-algebras. Let
$(Q, \ \searrow, \ \swarrow, \ \nearrow, \ \nwarrow)$ be a quadri-algebra. A 1-parameter formal deformation of $(Q, \ \searrow, \ \swarrow, \ \nearrow, \ \nwarrow)$ is the data of a quadri-algebra  $(Q[[h]], \ \searrow(h), \ \swarrow(h), \ \nearrow(h), \ \nwarrow(h))$ and a family of ennea-algebras $(T_\tau[[h]]:=(Q[[h]], \ \searrow^\tau(h), \ \swarrow^\tau(h), \ \nearrow^\tau(h), \ \nwarrow^\tau(h) \downarrow^\tau(h), \ \uparrow^\tau(h), \ \prec^\tau(h), \ \succ^\tau(h), \ \circ^\tau(h)))_{\tau \in k}$, such that,
$$\bar{\star}^\tau _{Ennea-algebra}(h) \longrightarrow \bar{\star}^\tau _{Quadri-algebra}(h) +  \uparrow^\tau(h)  +  \downarrow^\tau(h)  +   \circ^\tau(h) +  \succ^\tau(h) +  \prec^\tau(h),$$
with by definition:
$$x \ \bar{\star}^\tau _{Quadri-algebra}(h)  \ y:= x \nearrow^\tau(h) y+  x \searrow^\tau(h) y +  x \swarrow^\tau(h) y + x \nwarrow^\tau(h) y, \ \ \ \forall x,y \in Q[[h]],$$
which is associative for $\tau=0$.
 
The sums operations described in the definition of a $1$-ennea-algebra are obtained from the sums operations of a quadri-algebra (indicated  in parenthesis) by adding a complementary term:
\begin{eqnarray*}
\triangleright^\tau(h) &:=& (\nearrow^\tau(h)  +  \searrow^\tau(h) ) +  \succ^\tau(h)  , \\
\triangleleft^\tau(h)  &:=&  (\nwarrow^\tau(h)  +  \swarrow^\tau(h))  +  \prec^\tau(h), \\
\vee^\tau(h) &:=& (\searrow^\tau(h)  +  \swarrow^\tau(h))  +  \downarrow^\tau(h)  , \\
\wedge^\tau(h) &:=&  (\nearrow^\tau(h)  +  \nwarrow^\tau(h))  +    \uparrow^\tau(h)  , \\
\bar{\circ}^\tau(h) &:=&  \uparrow^\tau(h)  +  \downarrow^\tau(h)  +   \circ^\tau(h) .
\end{eqnarray*} 
For a given $\tau$, the additional operations $( \uparrow^\tau(h) , \ \ \downarrow^\tau(h),  \ \circ^\tau(h)))$ and $( \succ^\tau(h) , \  \prec^\tau(h), \  \circ^\tau(h))$ 
deforming the sums operations of a given quadri-algebra structure obey dendriform trialgebras axioms.
This 1-parameter deformation respects the symmetries, the vertical 
and the horizontal structures of the quadri-algebra and also
the structure of connected Hopf algebra of the augmented free quadri-algebra since this structure still holds for any $\tau \in k$.
The following picture summarizes our discussion. 
\begin{center}
\includegraphics*[width=8cm]{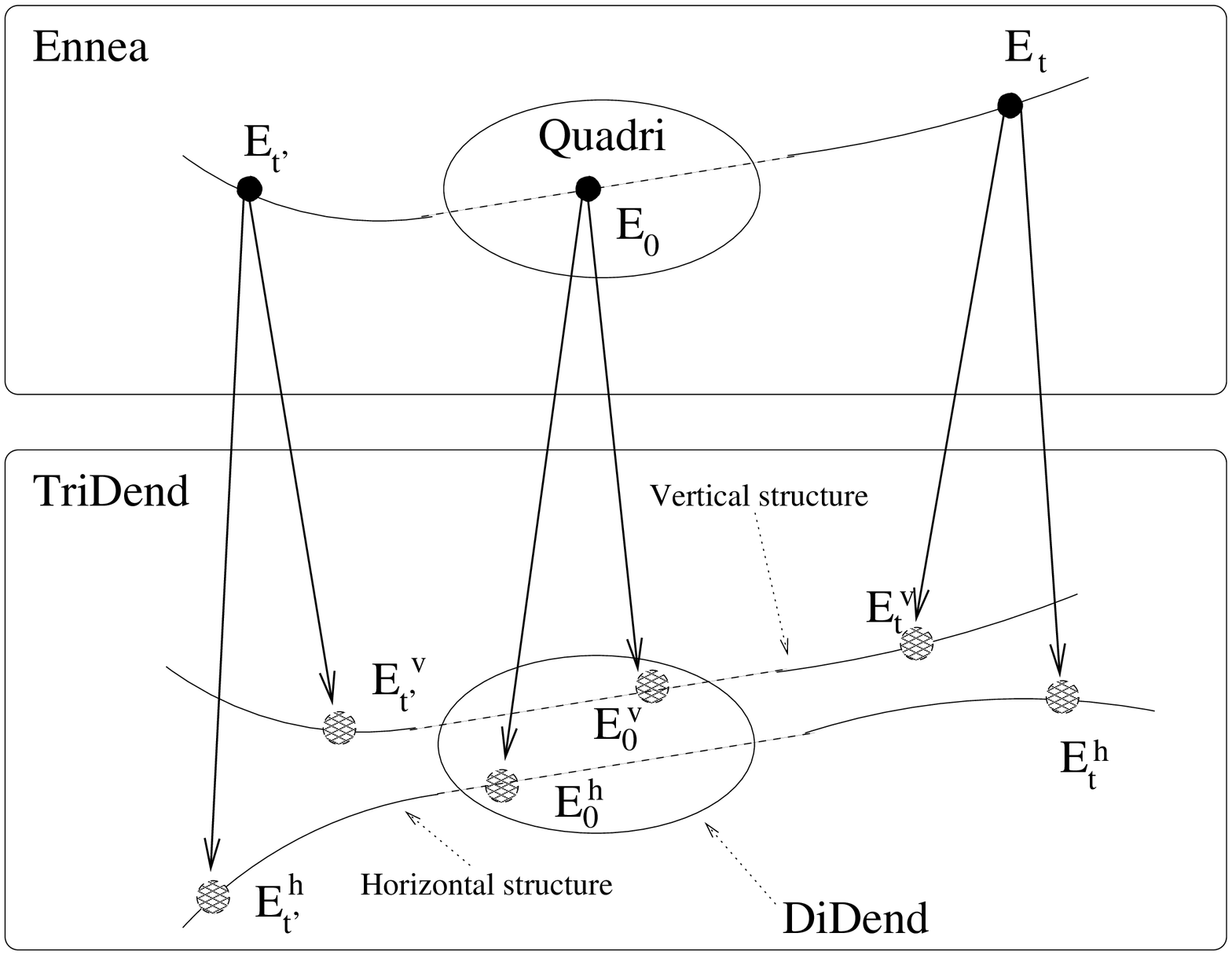}

\begin{scriptsize}
\textbf{
Family of ennea-algebras associated with a 1-parameter formal deformation of a quadri-algebra. Their vertical and horizontal structures are also represented.}
\end{scriptsize}
\end{center}
The same thing can be done in the category $\textsf{TriDend}$ of dendriform trialgebras. Let $(D, \ \prec, \ \succ)$ be a dendriform dialgebra. A 1-parameter deformation of $(D, \ \prec, \ \succ)$ is the data of a dendriform dialgebra $(D[[h]], \ \prec(h), \ \succ(h))$  and a family of dendriform trialgebras $(DT_\tau[[h]]:=(D[[h]], \ \prec^\tau(h), \ \succ^\tau(h), \ \circ^\tau(h))$ such that,
$$ \star_{Dendriform \ trialgebra}^\tau(h) \longrightarrow \star^\tau(h) _{Dendriform \ dialgebra}(h) +   \circ^\tau(h) ,$$
where for $\tau=0$, the operation $x \star^0 _{Dendriform \ dialgebra}(h) y := x \prec(h) y +  x \succ(h) y$ is associative.

In the following subsections, we will produce formal deformations of these special associative algebras, whose associative products are clusters of binary operations verifying dendriform di- or tri-algebra axioms, 
quadri-algebra axioms or ennea-algebra axioms. 
\subsection{Formal deformations of dendriform dialgebras}
Let $(D, \ \prec, \ \succ)$ be a dendriform dialgebra. Consider the dendriform dialgebra $(D[[h]], \  \tilde{\prec}, \ \tilde{\succ})$ where the operations $\tilde{\prec}, \ \tilde{\succ}: D[[h]] \otimes_{k[[h]]} D[[h]] \xrightarrow{} D[[h]]$ are defined for all $x:= \sum_{n \geq 1} x_n h^n$ and $y:= \sum_{n \geq 1} y_n h^n$ belongs to $D[[h]]$ by $x \tilde{\diamond} y := \sum_{n \geq 1} (\sum_{i+j:=n} \ x_i \diamond y_j)h^n$, where $\diamond \in \{ \prec, \ \succ \}$. Deform them by introducing three binary operations $\tilde{\prec}_1, \ \tilde{\succ}_1, \ \tilde{\circ}_1: D[[h]] \otimes_{k[[h]]} D[[h]] \xrightarrow{} D[[h]]$ as follows:
\begin{eqnarray*}
x \prec(h) y &:=& x \tilde{\prec} y + hx  \tilde{\prec}_1 y, \\
x \succ(h) y &:=& x \tilde{\succ} y + hx  \tilde{\succ}_1 y, \\
x \circ(h) y &:=& 0 + h x  \tilde{\circ}_1 y,
\end{eqnarray*}
for all $x,y \in D[[h]]$ and where $\tilde{\diamond}_1 \in \{ \tilde{\prec}_1, \ \tilde{\succ}_1, \ \tilde{\circ}_1 \}$ extends the operation $\diamond_1: D^{\otimes 2} \xrightarrow{} D$ to $D[[h]]$. 
Obviously, if $\tilde{\star} \longrightarrow \ \tilde{\succ} + \tilde{\prec}$, then set $\star(t) \longrightarrow \tilde{\star} + t\tilde{\star}_1$, where $\tilde{\star}_1 
\longrightarrow \ \tilde{\succ}_1 + \tilde{\prec}_1 + \tilde{\circ}_1$.
We require that $(D[[h]], \ \prec(h), \ \succ(h), \ \circ(h) )$ is a dendriform trialgebra.
To simplify the following equations, we suppose that  $\prec_1, \ \succ_1, \ \circ_1$ obey dendriform trialgebra axioms and
get,
\begin{eqnarray*}
(x \tilde{\prec}_1 y)\tilde{\prec} z + (x \tilde{\prec} y)\tilde{\prec}_1 z &=& x \tilde{\prec} ( y \tilde{\star}_1 z) + x \tilde{\prec}_1(y \tilde{\star} z),\\
(x \tilde{\succ}_1 y)\tilde{\prec} z  + (x \tilde{\succ} y)\tilde{\prec}_1 z &=& x \tilde{\succ}_1( y \tilde{\prec} z) + x \tilde{\succ}( y \tilde{\prec}_1 z),\\
(x \tilde{\star}_1 y)\tilde{\succ} z + (x \tilde{\star} y)\tilde{\succ}_1 z &=& x \tilde{\succ}_1(y \tilde{\succ} z) + x \tilde{\succ}(y \tilde{\succ}_1 z),\\
(x \tilde{\succ} y) \tilde{\circ}_1 z &=& x \tilde{\succ} (y \tilde{\circ}_1 z), \\
(x \tilde{\prec} y)\tilde{\circ}_1 z &=& x \tilde{\circ}_1(y \tilde{\succ} z),\\
(x \tilde{\circ}_1 y)\tilde{\prec} z &=& x \tilde{\circ}_1( y \tilde{\prec} z),
\end{eqnarray*}
for all $x,y,z \in D[[h]]$. This system of equations will be denoted by $\widetilde{(Syst. 1)}$. By $(Syst. 1)$, we just mean the same system of equations written without the notation tilde and thus with variables $x,y$ and $z$ belonging to $D$.
In general, $\widetilde{(Syst. 1)}$ is difficult to solve but can be addressed in the very important case where operations come from Baxter operators.
Suppose $\beta:D \xrightarrow{} D$ is a Baxter operator and $\beta_1: D \xrightarrow{} D$ is a $r_1$-Baxter operator with $r_1 \in k$. Then, for all $x,y,z \in D$, set as usual: $x \prec y := x\beta(y)$, $x \succ y := \beta(x)y$, $x \prec y := x\beta_1(y)$, $x \succ y := \beta_1(x)y$
and  $x \circ y := r_1xy$ and extend these operations to $D[[h]]$.
Then, our system of equations $\widetilde{(Syst. 1)}$ is reduced to the following equation:
$$ 
\widetilde{(Eq. \ 1)} \ \ \ \
\tilde{\beta}(x \tilde{\star}_1 y) + \tilde{\beta}_1(x \tilde{\star} y) := \tilde{\beta}_1(x)\tilde{\beta}(y) + \tilde{\beta}(x)\tilde{\beta}_1(y), \ \ \forall x,y \in D[[h]].
$$
That is:
$$ \widetilde{(Eq. \ 1)} \ \ \ \
\tilde{\beta}(x \tilde{\beta}_1(y)) + \tilde{\beta}( \tilde{\beta}_1(x)y) + r_1 \tilde{\beta}(xy) +  \tilde{\beta}_1(x  \tilde{\beta}(y)) +  \tilde{\beta}_1( \tilde{\beta}(x)y) :=  \tilde{\beta}_1(x) \tilde{\beta}(y) +  \tilde{\beta}(x) \tilde{\beta}_1(y), \ \ \forall x,y \in D[[h]],
$$
where Baxter operators $\beta$ and $\beta_1$ are respectively extended to $D[[h]] $ via the linear maps $\tilde{\beta}, \ \tilde{\beta}_1: D[[h]] \xrightarrow{} D[[h]]$ defined by $x:=\sum_{n \geq 0} x_n h^n \mapsto \tilde{\beta}(x):= \sum_{n \geq 0} \beta(x_n)h^n$ and $x:=\sum_{n \geq 0} x_n h^n \mapsto \tilde{\beta}_1(x):= \sum_{n \geq 0} \beta_1(x_n)h^n$. Observe that the linear map $\tilde{\beta}$ is a Baxter operator and $\tilde{\beta}_1$ is a $r_1$-Baxter operator. Moreover,
for all $x,y \in  D[[h]]$, we get $x \tilde{\prec} y := x\tilde{\beta}(y)$, $x \tilde{\succ} y := \tilde{\beta}(x)y$, $x \tilde{\prec}_1 y := x\tilde{\beta}_1(y)$, $x \tilde{\succ}_1 y := \tilde{\beta}_1(x)y$
and  $x \tilde{\circ}_1 y := r_1xy$.
\begin{theo}
\label{hyp. tri}
Let $(A, \ \mu, \ \Delta, \ \Delta_1)$ be a $k$-vector space such that:
\begin{enumerate}
\item {The $k$-vector space $(A, \ \mu, \ \Delta)$ is a $\epsilon$-bialgebra,}
\item {The $k$-vector space $(A, \ \mu, \ \Delta_1)$ is a $\epsilon(r_1)$-bialgebra,}
\item{$(\Delta_1 \otimes id ) \Delta = (id \otimes \Delta)\Delta_1$,}
\item{$(\Delta \otimes id ) \Delta_1 = (id \otimes \Delta_1)\Delta$.}
\end{enumerate} 
Then, the maps $\beta, \ \beta_1: (\textsf{End}(A), \ *, \ *_1) \xrightarrow{} (\textsf{End}(A), \ *, \ *_1)$, where $*$ and $*_1$ are the so-called convolution products associated with $\Delta$ and $\Delta_1$ respectively, defined by $\beta(T):= id * T$ and $\beta_1(T):= id *_1 T$, for any $T \in \textsf{End}(A)$
are respectively a Baxter operator and a $r_1$-Baxter operator. Moreover, 
$$(Eq. 1) \ \ \ 
\beta(T \beta_1(S)) + \beta(\beta_1(T)S) + r_1\beta(TS) + \beta_1(T \beta(S)) + \beta_1(\beta(T)S) := \beta_1(T)\beta(S) + \beta(T)\beta_1(S), 
$$
for all $T,S \in \textsf{End}(A)$.
\end{theo}
\Proof
Keep notation of Theorem \ref{hyp. tri}. The map $\beta$ is a Baxter operator since the $k$-vector space $(A, \ \mu, \ \Delta_1)$ is a $\epsilon$-bialgebra. 
Similarly, as the $k$-vector space $(A, \ \mu, \ \Delta_1)$ is a $\epsilon(r_1)$-bialgebra, the map $\beta_1$ is a $r_1$-Baxter operator. Fix
$T,S \in \textsf{End}(A)$ and
observe that Condition 3 implies,
$$ \beta_1(T)\beta(S) := \beta_1(T \beta(S)) + \beta(\beta_1(T)S) + r_1\beta(TS), $$
and that Condition 4 implies,
$$ \beta(T)\beta_1(S) := \beta(T \beta_1(S)) + \beta_1(\beta(T)S).$$
\eproof
\Rk
Theorem \ref{hyp. tri} also holds if we replace the operators $\beta$ and $\beta_1$ by $\gamma$ and $\gamma_1$ respectively.
\begin{exam}{[Weighted directed graphs]}
Such $k$-vector spaces verifying conditions of Theorem \ref{hyp. tri} exist. Let us give an example related to weighted directed graphs. Keep notation of Section \ref{examp}. From Proposition \ref{al}, $(kG, \ \Delta)$ is a $\epsilon$-bialgebra. From Proposition \ref{dup}, $(kG, \ \hat{\Delta})$ is a $\epsilon(-1)$-bialgebra. Now, observe that,
$$(\Delta \otimes id ) \hat{\Delta} = (id \otimes \hat{\Delta})\Delta, \ \ \ (\hat{\Delta} \otimes id ) \Delta = (id \otimes \Delta)\hat{\Delta}.$$
\end{exam}
\Rk
Such co-structures verifying,
 $$(\Delta_i \otimes id ) \Delta_j = (id \otimes \Delta_j)\Delta_i,$$
for all $i,j := 1,2$
are called $2$-hypercubic coalgebras in \cite{codialg1}, see also \cite{LodayRonco, Richter, BaxLer}. Such structures allow the construction of Hochschild 2-cocycles. Indeed, take the associated convolution products and observe that:
$$ (a *_2 b) *_1 c + (a *_1 b) *_2 c := a *_2(b *_1 c) + a *_1 (b *_2 c).$$ See \cite{BaxLer} for other examples of deformations linked to hypercubic structures.

Keep notation of Theorem \ref{hyp. tri}.
The operators $\beta$ and $\beta_1$ solutions of $(Eq. \ 1)$ give five operations on $\textsf{End}(A)$ defined for all $T, S \in \textsf{End}(A)$ by,
$ T \prec S := T \beta(S)$, $ T \succ S := \beta(T) S$ and $ T \prec_1 S := T \beta_1(S)$, $ T \succ_1 S := \beta_1(T) S$, $ T \circ_1 S := r_1T S$ which are solutions of $(Syst. \ 1)$. Extended to $\textsf{End}(A)[[h]]$, they give five binary operations defined for all $T, S \in \textsf{End}(A)[[t]]$ by,
$ T \tilde{\prec} S := T \tilde{\beta}(S)$, $ T \tilde{\succ} S := \tilde{\beta}(T) S$ and $ T \tilde{\prec}_1 S := T \tilde{\beta}_1(S)$, $ T \tilde{\succ}_1 S := \tilde{\beta}_1(T) S$, $ T \tilde{\circ}_1 S := r_1T S$ which are solutions of $\widetilde{(Syst. \ 1)}$. This formal deformation is denoted by $TD_1[[h]]:=(\textsf{End}(A)[[h]], \prec^1(h), \succ^1(h), \circ^1(h))$ where,
\begin{eqnarray*}
T \prec^1(h) S &:=& T \tilde{\prec} S + hT  \tilde{\prec}_1 S, \\
T \succ^1(h) S &:=& T \tilde{\succ} S + hT  \tilde{\succ}_1 S, \\
T \circ^1(h) S &:=& 0 + h T  \tilde{\circ}_1 S,
\end{eqnarray*}
 
Proposition \ref{tensor} turns $(\textsf{End}(A), \ \prec_{\beta}, \ \succ_{\beta})^{ \otimes  2}$ into a quadri-algebra and $(TD_1[[h]]:=(D[[h]], \prec^1(h), \succ^1(h), \circ^1(h)))^{ \otimes_{k[[h]]}  2} $ into a $1$-ennea-algebra denoted by $ETD_1[[h]]$. (Proposition \ref{tensor} still holds if $\otimes$ is replaced by $\otimes_{k[[h]]}$ since if $A$ and $B$ are two $k$-vector spaces, $(A \otimes B)[[h]]=A[[h]] \otimes_{k[[h]]} B[[h]]$.)
By applying  Proposition \ref{Scale} to the $r_1$-Baxter operator $\beta_1$, a family $(\beta^\tau:= \tau\beta_1)$ of $\tau r_1$-Baxter operators is obtained as well as a family of ennea-algebras $ETD_\tau[[h]]:=(D[[h]], \prec^\tau(h), \succ^\tau(h), \circ^\tau(h))^{ \otimes_{k[[h]]}  2} $. For $\tau =0$, all the operations labelled by 1 vanish and $ETD_0[[h]]:=(D[[h]], \prec^0(h), \succ^0(h), \circ^0(h))^{ \otimes_{k[[h]]}  2} $ is isomorphic, as quadri-algebra, to the trivial formal deformation of $(\textsf{End}(A), \ \prec_{\beta}, \ \succ_{\beta})^{ \otimes  2}$.
Let us give an example of operations associated with $ETD_\tau[[h]]$. Via Proposition \ref{tensor}, we get:
\begin{eqnarray*}
(a \otimes_{k[[h]]} b) \searrow^\tau(h) (a' \otimes_{k[[h]]} b') &:=& (a \otimes_{k[[h]]} b) \searrow^\tau (a' \otimes_{k[[h]]} b') + h (a \otimes_{k[[h]]} b) \searrow^\tau _1 (a' \otimes_{k[[h]]} b') \\ & & + h^2 (a \otimes_{k[[h]]} b) \searrow^\tau _2 (a' \otimes_{k[[h]]} b'),
\end{eqnarray*}
where,
$$ (a \otimes_{k[[h]]} b) \searrow^\tau (a' \otimes_{k[[h]]} b'):= (a \prec a') \otimes_{k[[h]]} (b \prec b'), \ \ \
(a \otimes_{k[[h]]} b) \searrow^\tau _2 (a' \otimes_{k[[h]]} b'):= (a \prec^\tau _1 a') \otimes_{k[[h]]} (b \prec^\tau_1 b'),$$
$$(a \otimes_{k[[h]]} b) \searrow^\tau _1 (a' \otimes_{k[[h]]} b'):= (a \prec a') \otimes_{k[[h]]} (b \prec^\tau _1 b') + (a \prec^\tau _1 a') \otimes_{k[[h]]} (b \prec b'),$$
with $a,a',b,b' \in \textsf{End}(A)[[h]]$. Recall for instance that $b \prec^\tau_1 b' := \tau b \beta(b')$.
\Rk
A binary, quadratic and non-symmetric operad can be associated with this formal deformation. It has 5 binary operations which are $\prec, \ \succ, \ \prec_1, \ \succ_1, \ \circ_1$ and has 16 relations. Indeed, 3 relations are needed for $\prec, \ \succ$ verify dendriform dialgebras axioms, 7 are needed for $\prec_1, \ \succ_1, \ \circ_1$ verify dendriform trialgebras axioms and 6 relations describing $(Syst. \ 1)$ are required to establish the formal deformation. The sequence of the dimensions associated with this operad starts with $1,5,34$.
\Rk 
The case $r_1 := 0$ is easily obtained. A binary, quadratic and non-symmetric operad can also be associated with this formal deformation. It has 4 binary operations which are $\prec, \ \succ, \ \prec_1, \ \succ_1$ and has 9 relations. Indeed, 3 relations are needed for $\prec, \ \succ$ verify dendriform dialgebras axioms, as well as 3 for $\prec_1, \ \succ_1$. However, $(Syst. \ 1)$ is reduced to the first 3 relations. The sequence of the dimensions associated with this operad starts with $1,4,23$. 
\Rk
There exists a commutative version of these two types of algebras. Indeed, the relations defining these two operads are globally invariant
under the following transformation:
$x \prec^{op} y := y \succ x$, $x \succ^{op} y := y \prec x$, and $x \prec_1 ^{op} y := y \succ_1 x$, $x \succ_1 ^{op} y := y \prec_1 x$, $x \circ_1 ^{op} y := y \circ_1 x$. These algebras are said to be {\it{commutative}} when they coincide with their opposite structures.
\begin{prop}
\label{zz}
With the following choice of the unit action: $x \prec 1 = x$ and $1 \succ x = x$, (and all the other actions equal to zero) the augmented free algebras associated with the two operads just defined, as well as their commutative versions, have a structure of connected Hopf algebra. 
\end{prop}
\Proof
It suffices to prove that this choice is compatible and coherent with the axioms defining the operads just described, which is straightforward. As regards their commutative versions, observe that $x \prec^{op} 1 := 1 \succ x := x$ and $1 \succ^{op} x := x \prec 1:=x$, the choice of the unit action is in agreement with the opposite structure.
\eproof
\begin{prop}
\label{zz1}
With the following choice of the unit action: $x \prec_1 1 = x$ and $1 \succ_1 x = x$, (and all the other actions equal to zero) the augmented free algebras, as well as their commutative versions, associated with the two operads just defined have a structure of connected Hopf algebra. 
\end{prop}
\Proof
Observe that the non-null equations are symmetric when the operations labelled by 1 are replaced by the same operations with no labels and conversely.
\eproof

\noindent
Observe that the operation $\bar{\star} \longrightarrow \star + \star_1$
is associative. In Proposition \ref{zz}, $1 \bar{\star} x=x= x \bar{\star} 1$ and $1 \star x=x= x \star 1$. In Proposition \ref{zz1}, $1 \bar{\star} x=x= x \bar{\star} 1$ and $1 \star_1 x=x= x \star_1 1$. Therefore,
in the equations $(1)$ and $(2)$ p.14, we can choose either the associative product  $\bar{\star}$ or $\star$ for the Hopf algebra structure of Proposition \ref{zz}. We can choose either the associative product  $\bar{\star}$ or $\star_1$ for the Hopf algebra structure of Proposition \ref{zz1}. This gives at least four connected Hopf algebras on the augmented free algebras associated with the operads described in these two Remarks.\begin{exam}{[Constructions of $\epsilon$-bialgebras]}
\begin{prop}
Let $V$ and $C$ be two sets. Consider the free $k$-vector space generated by $V$, denoted by $kV$. Suppose $V$ is equipped with a set $Co$ of co-operations verifying a set $S$ of relations of type $(\Delta_j \otimes id)\Delta_i :=(id \otimes \Delta_i)\Delta_j$, with $\Delta_i, \ \Delta_j: kV \xrightarrow{} kV^{\otimes 2}$, for all $\Delta_i, \ \Delta_j \in Co$. Denote by $As(C)$ the free associative algebra generated by $C$. Let $\Delta_i \in Co$. Extend $\Delta_i$ as follows: $\Delta_i(a)=0$, for all $a \in As(C)$ and $\Delta_i(ab)=\Delta_i(a)b + a\Delta_i(b)$, for all $a,b \in As(C)\bra V \ket$, the non-commutative polynomial algebra with indeterminates belonging to $V$ and coefficients belonging to $As(C)$. Then, the set $S$ of relations still holds in $As(C)\bra V \ket$.
\end{prop}
\Proof
Straightforward.
\eproof

\noindent
Apply this proposition to tilings of the $(n^2,1)$-De Bruijn graphs, $n>1$, found in \cite{codialg1} to produce explicit examples of $k$-vector spaces verifying conditions of Theorem \ref{hyp. tri} with $r_1 := 0$.
\begin{center}
\includegraphics*[width=3cm]{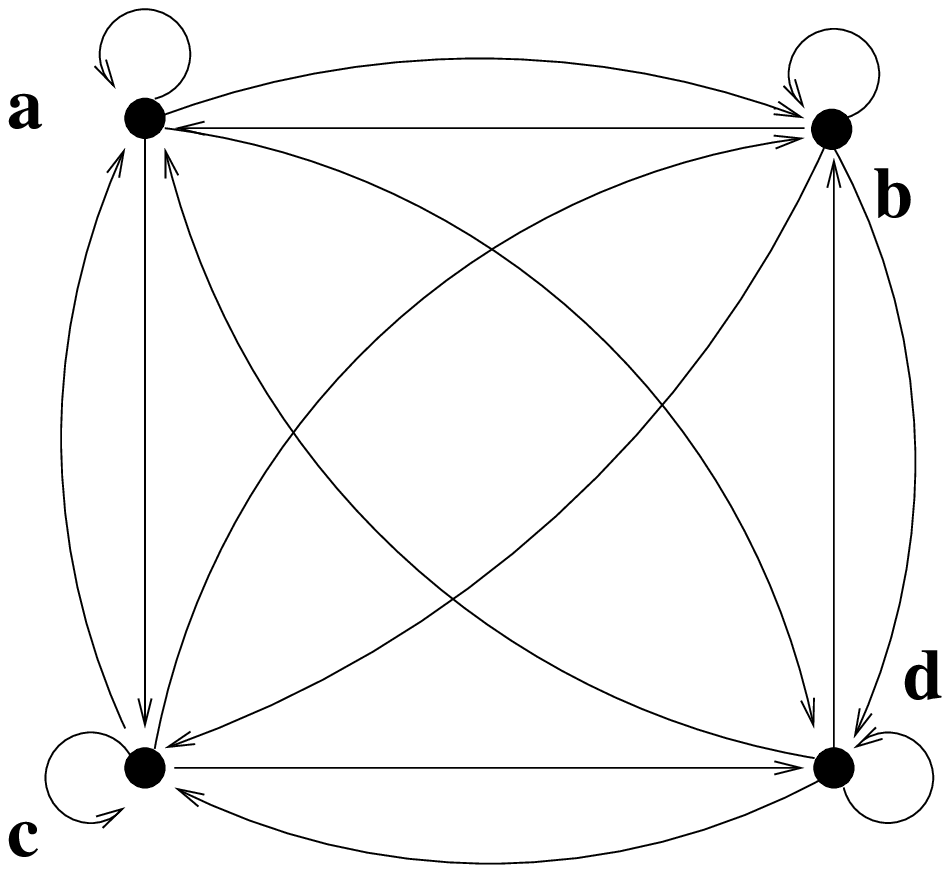}

\begin{scriptsize} \textbf{The $(4,1)$-De Bruijn graph. }
\end{scriptsize}
\end{center}

\end{exam}
\begin{exam}{}
Other examples can be found in the work of Voiculescu \cite{Voi} on free probability. 
\end{exam}

\subsection{Formal deformations of dendriform trialgebras}
To deform a dendriform trialgebra $(T, \ \prec, \ \succ, \ \circ )$, three other operations, $\prec_1, \ \succ_1, \ \circ_1: T^{\otimes 2} \xrightarrow{} T$, verifying dendriform trialgebra axioms are needed. Therefore, we are looking for solving the system of 7 equations, 
$$x \diamond(h) y:= x \tilde{\diamond} y +  hx \tilde{\diamond}_1 y,$$
for all $x,y \in T[[h]]$ and
where $\tilde{\diamond}, \ \tilde{\diamond}_1: T[[h]] \otimes_{k[[h]]} T[[h]] \xrightarrow{} T[[h]]$ extend respectively the operations $\diamond$ and $\diamond_1$, with $\diamond \in \{\prec, \ \succ, \ \circ \}$ and $\diamond_1 \in \{\prec_1, \ \succ_1, \ \circ_1 \}$. Requiring that $ \prec(h), \ \succ(h)$ and $\circ(h) $ obey dendriform trialgebra axioms implies the system of equations $\widetilde{(Syst. \ 2)}$,
\begin{eqnarray*}
(x \tilde{\prec}_1 y)\tilde{\prec} z + (x \tilde{\prec} y)\tilde{\prec}_1 z &=& x \tilde{\prec} ( y \tilde{\star}_1 z) + x \tilde{\prec}_1(y \tilde{\star} z),\\
(x \tilde{\succ}_1 y)\tilde{\prec} z  + (x \tilde{\succ} y)\tilde{\prec}_1 z &=& x \tilde{\succ}_1( y \tilde{\prec} z) + x \tilde{\succ}( y \tilde{\prec}_1 z),\\
(x \tilde{\star}_1 y)\tilde{\succ} z + (x \tilde{\star} y)\tilde{\succ}_1 z &=& x \tilde{\succ}_1(y \tilde{\succ} z) + x \tilde{\succ}(y \tilde{\succ}_1 z),\\
(x \tilde{\succ} y) \tilde{\circ}_1 z  + (x \tilde{\succ}_1 y) \tilde{\circ} z &=& x \tilde{\succ} (y \tilde{\circ}_1 z) + x \tilde{\succ}_1 (y \tilde{\circ} z), \\
(x \tilde{\prec} y)\tilde{\circ}_1 z + (x \tilde{\prec}_1 y)\tilde{\circ} z&=& x \tilde{\circ}_1(y \tilde{\succ} z) + x \tilde{\circ}(y \tilde{\succ}_1 z),\\
(x \tilde{\circ}_1 y)\tilde{\prec} z + (x \tilde{\circ} y)\tilde{\prec}_1 z &=& x \tilde{\circ}_1( y \tilde{\prec} z) + x \tilde{\circ}( y \tilde{\prec}_1 z), \\
(x \tilde{\circ}_1 y)\tilde{\circ} z + (x \tilde{\circ} y)\tilde{\circ}_1 z &=& x \tilde{\circ}_1( y \tilde{\circ} z) + x \tilde{\circ}( y \tilde{\circ}_1 z),
\end{eqnarray*}
for all $x,y,z \in T[[h]]$.
These equations can be solved when operations $\prec, \ \succ, \ \circ$ comes from a $r$-Baxter operator and $\prec_1, \ \succ_1, \ \circ_1$ from a $r_1$-Baxter operator. 
Suppose $\beta:T \xrightarrow{} T$ is a $r$-Baxter operator and $\beta_1: T \xrightarrow{} T$ is a $r_1$-Baxter operator, $r,r_1 \in k$ different from zero. Then, for all $x,y,z \in T$, set: $x \prec y := x\beta(y)$, $x \succ y := \beta(x)y$, $x \circ y := rxy$ and $x \prec_1 y := x\beta_1(y)$, $x \succ_1 y := \beta_1(x)y$, $x \circ_1 y := r_1xy$. Extend these operations to $T[[h]]$.
Then, our system of equations $\widetilde{(Syst. \ 2)}$ is reduced to the following equation:
$$ 
\widetilde{(Eq. \ 2)} \ \ \ \
\tilde{\beta}(x \tilde{\star}_1 y) + \tilde{\beta}_1(x \tilde{\star} y) := \tilde{\beta}_1(x)\tilde{\beta}(y) + \tilde{\beta}(x)\tilde{\beta}_1(y), \ \ \forall x,y \in T[[h]].
$$
That is:
$$ \widetilde{(Eq. \ 2)} \ \ \ \
\tilde{\beta}(x \tilde{\beta}_1(y)) + \tilde{\beta}( \tilde{\beta}_1(x)y) + r_1 \tilde{\beta}(xy) + r \tilde{\beta_1}(xy) + \tilde{\beta}_1(x  \tilde{\beta}(y)) +  \tilde{\beta}_1( \tilde{\beta}(x)y) :=  \tilde{\beta}_1(x) \tilde{\beta}(y) +  \tilde{\beta}(x) \tilde{\beta}_1(y), 
$$
for all $x,y \in T[[h]]$,
where Baxter operators $\beta$ and $\beta_1$ are respectively extended to $T[[h]] $ via the linear maps $\tilde{\beta}, \ \tilde{\beta}_1: T[[h]] \xrightarrow{} T[[h]]$ defined by $x:=\sum_{n \geq 0} x_n h^n \mapsto \tilde{\beta}(x):= \sum_{n \geq 0} \beta(x_n)h^n$ and $x:=\sum_{n \geq 0} x_n h^n \mapsto \tilde{\beta}_1(x):= \sum_{n \geq 0} \beta_1(x_n)h^n$. Observe that the linear map $\tilde{\beta}$ is a $r$-Baxter operator and $\tilde{\beta}_1$ is a $r_1$-Baxter operator. Moreover,
for all $x,y \in  T[[h]]$, we get $x \tilde{\prec} y := x\tilde{\beta}(y)$, $x \tilde{\succ} y := \tilde{\beta}(x)y$, $x \tilde{\circ} y := rxy$ and $x \tilde{\prec}_1 y := x\tilde{\beta}_1(y)$, $x \tilde{\succ}_1 y := \tilde{\beta}_1(x)y$
and  $x \tilde{\circ}_1 y := r_1 xy$.
Solving $\widetilde{(Syst. \ 2)}$ in this case leads to an analogue of Theorem \ref{hyp. tri}.
\begin{theo}
\label{hyp. tri2}
Let $(A, \ \mu, \ \Delta, \ \Delta_1)$ be a $k$-vector space such that:
\begin{enumerate}
\item {The $k$-vector space $(A, \ \mu, \ \Delta)$ is a $\epsilon(r)$-bialgebra, $r \in k$ different from zero,}
\item {The $k$-vector space $(A, \ \mu, \ \Delta_1)$ is a $\epsilon(r_1)$-bialgebra, $r_1 \in k$ different from zero,}
\item{$(\Delta_1 \otimes id ) \Delta = (id \otimes \Delta)\Delta_1$,}
\item{$(\Delta \otimes id ) \Delta_1 = (id \otimes \Delta_1)\Delta$.}
\end{enumerate} 
Then, $\beta, \ \beta_1: (\textsf{End}(A), \ *, \ *_1) \xrightarrow{} (\textsf{End}(A), \ *, \ *_1)$, where $*$ and $*_1$ are the so-called convolution products associated with $\Delta$ and $\Delta_1$ respectively, defined by $\beta(T):= id * T$ and $\beta_1(T):= id *_1 T$, for any $T \in \textsf{End}(A)$
are respectively a $r$-Baxter operator and a $r_1$-Baxter operator. Moreover, we get $(Eq. \ 2)$,
$$ 
\beta(T \beta_1(S)) + \beta(\beta_1(T)S) + r_1\beta(TS) + r \beta_1(TS) + \beta_1(T \beta(S)) + \beta_1(\beta(T)S) := \beta_1(T)\beta(S) + \beta(T)\beta_1(S), 
$$
for all $T,S \in \textsf{End}(A)$.
\end{theo}
\Proof
Similar to Theorem \ref{hyp. tri}.
\eproof

\noindent
Keep notation of Theorem \ref{hyp. tri2}.
The operators $\beta$ and $\beta_1$ solutions of $(Eq. \ 2)$ give six operations on $\textsf{End}(A)$ defined for all $T, S \in \textsf{End}(A)$ by,
$ T \prec S := T \beta(S)$, $ T \succ S := \beta(T) S$ and $ T \prec_1 S := T \beta_1(S)$, $ T \succ_1 S := \beta_1(T) S$, $ T \circ_1 S := r_1T S$. Therefore, these operations are solutions of $(Syst. \ 2)$. Extended to $\textsf{End}(A)[[h]]$, they give five operations defined for all $T, S \in \textsf{End}(A)[[h]]$ by,
$ T \tilde{\prec} S := T \tilde{\beta}(S)$, $ T \tilde{\succ} S := \tilde{\beta}(T) S$ and $ T \tilde{\prec}_1 S := T \tilde{\beta}_1(S)$, $ T \tilde{\succ}_1 S := \tilde{\beta}_1(T) S$, $ T \tilde{\circ}_1 S := r_1T S$ which are solutions of $\widetilde{(Syst. \ 2)}$. This formal deformation is denoted by $TrD_1[[h]]:=(\textsf{End}(A)[[h]], \prec^1(h), \succ^1(h), \circ^1(h))$ where,
\begin{eqnarray*}
T \prec^1(h) S &:=& T \tilde{\prec} S + hT  \tilde{\prec}_1 S, \\
T \succ^1(h) S &:=& T \tilde{\succ} S + hT  \tilde{\succ}_1 S, \\
T \circ^1(h) S &:=& T  \tilde{\circ} S + h T  \tilde{\circ}_1 S.
\end{eqnarray*}
By applying Proposition \ref{Scale} on the opearator $\beta_1$, a family
of dendriform trialgebras $(TrD_\tau [[h]]:=(\textsf{End}(A)[[h]], \prec^\tau(h), \succ^\tau(h), \circ^\tau(h)))_{\tau \in k}$ is created. Moreover, $TrD_0 [[h]]:=(\textsf{End}(A)[[h]], \prec^0(h), \succ^0(h), \circ^0(h))$ is isomorphic, as dendriform trialgebras, to the trivial deformation of $(\textsf{End}(A), \prec, \succ, \circ)$.
\Rk
A binary, quadratic and non-symmetric operad can be associated with this deformation. The sequence of its dimensions starts with $1,6,51$.
Here again, the choice of the unit action defined in Proposition \ref{zz}
or \ref{zz1} turns the augmented free algebra associated with the operad just described into a connected Hopf algebra.
There also exists a commutative version of these two types of algebras since the relations defining them are globally invariant
under the following transformation:
$x \prec^{op} y := y \succ x$, $x \succ^{op} y := y \prec x$, $x \circ^{op} y := y \circ x$ and $x \prec_1 ^{op} y := y \succ_1 x$, $x \succ_1 ^{op} y := y \prec_1 x$, $x \circ_1 ^{op} y := y \circ_1 x$. These algebras are said to be {\it{commutative}} when they coincide with their opposite structures. The augmented free commutative algebra associated with the operad just described has also a structure of connected Hopf algebra.
\subsection{Formal deformations of quadri-algebras}
Let $(Q, \ \nearrow, \ \searrow, \ \swarrow, \ \nwarrow)$ be a quadri-algebra. Extend as usual the operations $\diamond : Q^{\otimes 2} \xrightarrow{} Q$ to $\tilde{\diamond} : Q[[h]] \otimes_{k[[h]]} Q[[h]] \xrightarrow{} Q[[h]]$, for all $\diamond:=\ \nearrow, \ \searrow, \ \swarrow, \ \nwarrow$. Consider four other operations $\nearrow_1, \ \searrow_1, \ \swarrow_1, \ \nwarrow_1: Q^{\otimes 2} \xrightarrow{} Q$ extended to $\tilde{\nearrow}_1, \ \tilde{\searrow}_1, \ \tilde{\swarrow}_1, \ \tilde{ \nwarrow}_1: Q[[h]] \otimes_{k[[h]]} Q[[h]] \xrightarrow{} Q[[h]]$ and set:
$$ x \nearrow(h) y := x \tilde{\nearrow} y + hx \tilde{\nearrow}_1 y ,$$
$$ x \searrow(h) y := x \tilde{\searrow} y + hx \tilde{\searrow}_1 y ,$$
$$ x \swarrow(h) y := x \tilde{\swarrow} y + hx \tilde{\swarrow}_1 y ,$$
$$ x \nwarrow(h) y := x \tilde{\nwarrow} y + hx \tilde{\nwarrow}_1 y .$$
Suppose $\nearrow_1, \ \searrow_1, \ \swarrow_1, \ \nwarrow_1$ obey quadri-algebras axioms. The sums operations can be also written, for instance $\bar{\star}(h) \longrightarrow \tilde{\bar{\star}} + h\tilde{\bar{\star}}_1,$ and so on.
We require that $(Q[[h]], \ \nearrow(h), \ \searrow(h), \ \swarrow(h), \ \nwarrow(h))$ is a quadri-algebra. This implies that operations labelled by 1 and operations with no label have to verify the following system of equations $\widetilde{(Syst. \ 3)}$:
\begin{eqnarray*}
(x \tilde{\nwarrow}_1 y)\tilde{\nwarrow} z + (x \tilde{\nwarrow} y)\tilde{\nwarrow}_1 z &=& x \tilde{\nwarrow}(y \tilde{\bar{\star}}_1 z) + x \nwarrow_1(y \tilde{\bar{\star}} z), \\
(x \nearrow y) \tilde{\nwarrow}_1 z + (x \tilde{\nearrow}_1 y) \tilde{\nwarrow} z &=& x\tilde{\nearrow}(y \tilde{\triangleleft}_1 z) + x \tilde{\nearrow}_1(y \tilde{\triangleleft} z), \\
(x \tilde{\wedge}_1 y)\tilde{\nearrow}  z+ (x \tilde{\wedge} y)\tilde{\nearrow}_1 z &=& x \tilde{\nearrow}(y \tilde{\triangleright}_1 z) + x \tilde{\nearrow}_1(y \tilde{\triangleright} z), \\
(x \tilde{\swarrow}_1 y)\tilde{\nwarrow} z + (x \tilde{\swarrow} y)\tilde{\nwarrow}_1 z &=& x \tilde{\swarrow}(y \tilde{\wedge}_1 z) +x \swarrow_1(y \tilde{\wedge} z),\\
(x \tilde{\searrow}_1 y )\tilde{\nwarrow} z +( x \tilde{\searrow} y)\tilde{\nwarrow}_1 z &=& x \tilde{\searrow}_1(y \tilde{\nwarrow} z) +x \tilde{\searrow}(y \tilde{\nwarrow}_1 z),\\
(x \tilde{\vee}_1 y)\tilde{\nearrow} z +  (x \tilde{\vee} y)\tilde{\nearrow}_1 z &=& x \tilde{\searrow}(y \tilde{\nearrow}_1 z) + x \tilde{\searrow}_1(y \tilde{\nearrow} z),\\
(x \tilde{\triangleleft} y)\tilde{\swarrow}_1 z + (x \tilde{\triangleleft}_1 y)\tilde{\swarrow} z &=& x \tilde{\swarrow}_1(y \tilde{\vee} z) + x \tilde{\swarrow}(y \tilde{\vee}_1 z),\\
(x \tilde{\triangleright} y)\tilde{\swarrow}_1 z + (x \tilde{\triangleright}_1 y)\tilde{\swarrow }z &=& x \tilde{\searrow}(y \tilde{\swarrow}_1 z) + x \tilde{ \searrow}_1(y \tilde{\swarrow} z),\\
(x \tilde{\bar{\star}} y)\tilde{\searrow}_1 z + (x \tilde{\bar{\star}}_1 y)\tilde{\searrow} z &=& x \tilde{\searrow}(y \tilde{\searrow}_1 z)
+ x \tilde{\searrow}_1(y \tilde{\searrow} z).
\end{eqnarray*}
To solve this system in the particular case where operations come from Baxter operators, the following theorem is needed.
\begin{theo}
\label{hyp. tri3}
Let $(A, \ \mu, \ \Delta, \ \Delta_1)$ be a $k$-vector space such that:
\begin{enumerate}
\item {The $k$-vector space $(A, \ \mu, \ \Delta)$ and $(A, \ \mu, \ \Delta_1)$ are both $\epsilon$-bialgebras,}
\item{We suppose $(\Delta_1 \otimes id ) \Delta = (id \otimes \Delta)\Delta_1$ and $(\Delta \otimes id ) \Delta_1 = (id \otimes \Delta_1)\Delta$.}
\end{enumerate} 
Then, $\beta, \gamma, \beta_1: (\textsf{End}(A), \ *, \ *_1) \xrightarrow{} (\textsf{End}(A), \ *, \ *_1)$, where $*$ and $*_1$ are the so-called convolution products associated with $\Delta$ and $\Delta_1$ respectively, defined by $\beta(T):= id * T$, $\gamma(T):=  T * id$  and $\beta_1(T):= id *_1 T$, for any $T \in \textsf{End}(A)$,
are Baxter operators which obey commutation relations $\beta\gamma :=\gamma\beta $ and $\beta_1\gamma :=\gamma\beta_1 $. Introduce the following operations:
$$T \searrow S := \gamma(T) \succ_{\beta} S, \ \ T \nearrow S:= T \succ_{\beta} \gamma(S), \ \ T \swarrow S := \gamma(T) \prec_{\beta} S, \ \ T \nwarrow S := T \prec_{\beta} \gamma(S),$$
and,
$$T \searrow_1 S := \gamma(T) \succ_{\beta_1} S, \ \ T \nearrow_1 S:= T \succ_{\beta_1} \gamma(S), \ \ T \swarrow_1 S := \gamma(T) \prec_{\beta_1} S, \ \ T \nwarrow_1 S := T \prec_{\beta_1} \gamma(S).$$
Then, with these operations the system of equations $(Syst. \ 3)$ can be solved.  
\end{theo}
\Proof
Let $(A, \ \mu, \ \Delta, \ \Delta_1)$ be a $k$-vector space verifying conditions of Theorem \ref{hyp. tri3}. Via the two Baxter operators $\beta$ and $\beta_1$, two dendriforms dialgebras are introduced on $\textsf{End}(A)$ . The binary operations are then indicated by $\prec_{\beta}, \ \succ_{\beta}$ on the one hand and by 
$\prec_{\beta_1}, \ \succ_{\beta_1}$, on the other hand.
The operator $\beta$ commutes with $\gamma$ which generates a quadri-algebra $(A, \ \searrow, \ \nearrow, \ \swarrow, \ \nwarrow) $ whose operations are defined in Theorem \ref{hyp. tri3}. Similarly, $\beta_1$ commutes with $\gamma$ which generates a quadri-algebra $(\textsf{End}(A), \ \searrow_1, \ \nearrow_1, \ \swarrow_1, \ \nwarrow_1) $ whose operations are defined in Theorem \ref{hyp. tri3}. Now, observe that the system of equations $Syst. 3'$
\begin{eqnarray*}
(R \prec_{\beta_1} S)\prec_{\beta} T + (R \prec_{\beta} S)\prec_{\beta_1} T &= & R \prec_{\beta} ( S \star_{\beta_1} T) + R \prec_{\beta_1}(S \star_{\beta} T), \\
(R \succ_{\beta_1} S)\prec_{\beta} T  + (R \succ_{\beta} S)\prec_{\beta_1} T &=& R \succ_{\beta_1}( S \prec_{\beta} T) + R \succ_{\beta}( S \prec_{\beta_1} T),\\
(R \star_{\beta_1} S)\succ_{\beta} T + (R \star_{\beta} S)\succ_{\beta_1} T &=& R \succ_{\beta_1}(S \succ_{\beta} T) + R \succ_{\beta}(S \succ_{\beta_1} T),
\end{eqnarray*}
holds for all $R,S,T \in \textsf{End}(A)$ thanks to Theorem \ref{hyp. tri}, with $r_1 =0$. 
Consider the first equation and replace $S$ by $\gamma(S)$ and $T$ by $\gamma(T)$ we obtain:
$$
(R \prec_{\beta_1} \gamma(S))\prec_{\beta} \gamma(T) + (R \prec_{\beta} \gamma(S))\prec_{\beta_1} \gamma(T) = R \prec_{\beta} ( \gamma(S) \star_{\beta_1} \gamma(T)) + R \prec_{\beta_1}(\gamma(S) \star_{\beta} \gamma(T)),
$$
which is equal to:
$$
(R \prec_{\beta_1} \gamma(S))\prec_{\beta} \gamma(T) + (R \prec_{\beta} \gamma(S))\prec_{\beta_1} \gamma(T) = R \prec_{\beta}  \gamma(S \star_{\beta_1} \gamma(T) + \gamma(S)\star_{\beta_1} T)  + R \prec_{\beta_1}\gamma(S \star_{\beta} \gamma(T) + \gamma(S)\star_{\beta} T),
$$
which is equal to:
$$(R \nwarrow_1 S)\nwarrow T + (R \nwarrow S)\nwarrow_1 T = R \nwarrow(S \bar{\star}_1 T) + R \nwarrow_1(S \bar{\star} T).$$
We obtain three equations from a given one of $(Syst. \ 3')$ by playing on the replacements of two variables, for intance $R$ by $\gamma(R)$ and $S$ by $\gamma(S)$ and so on. The nine equations of the system $(Syst. \ 3)$ are then obtained. Once our six operations are extended to $\textsf{End}(A)[[h]]$, the system $\widetilde{(Syst. \ 3)}$ will be solved.
\eproof

\noindent
By applying Proposition \ref{Scale} to the operator $\beta_1$, a family of quadri-algebras $(Q_\tau[[h]]:=(\textsf{End}(A)[[h]], \ \nearrow^\tau(h), \ \searrow^\tau(h), \ \swarrow^\tau(h), \ \nwarrow^\tau(h))_{\tau \in k}$ is obtained. When $\tau =0$, $\beta_1 ^{\tau}:= \tau \beta_1=0$ giving the quadri-algebra $Q_{0}[[t]]$ which is the trivial deformation of $Q$. 
\Rk
Another formal deformation from Theorem \ref{hyp. tri3} can be obtained by observing that the Baxter operator $\gamma$ and $\gamma_1$ commute with $\beta$.
\Rk
A binary, quadratic and non-symmetric operad can be obtained. The sequence of these dimensions starts with $1,8,101$. Such a start is also obtained by studying octo-algebras \cite{BaxLer}. These associative algebras have associative products which can be decomposed in four operations verifying quadri-algebras axioms in three different ways. 
By choosing the action of the unit as follows, $x \nwarrow 1 := x, 1 \searrow x := x$ or $x \nwarrow_1 1 := x, 1 \searrow_1 x := x$ (and all the other actions equal to zero) the augmented free algebra associated with this operad has a connected Hopf algebra structure. Observe also that there exists a notion of transpose and opposite associated with this type of algebras. The transformations are the same that those describing the transpose and opposite of a quadri-algebra and have to be applied to the unlabelled operations as well as to those labelled by 1. As our choice of the unit action is in agreement with the opposite structure,
the augmented free commutative algebra will have a structure of connected Hopf algebra.
\subsection{Formal deformations of ennea-algebras}
Let $(E, \ \nearrow, \ \searrow, \ \swarrow, \ \nwarrow, \ \downarrow, \ \uparrow, \ \prec, \ \succ, \ \circ)$ be a $r$-ennea-algebra. Denote by $I_E := \{ \nearrow, \ \searrow, \ \swarrow, \ \nwarrow, \ \nwarrow, \ \downarrow, \ \uparrow, \ \prec, \ \succ, \ \circ \}$.
Extend as usual the operations $\diamond : E^{\otimes 2} \xrightarrow{} E$ to $\tilde{\diamond} : E[[h]] \otimes_{k[[h]]} E[[h]] \xrightarrow{} E[[h]]$, for all $\diamond \in I_E$. Consider also nine other operations $\nearrow_1, \ \searrow_1, \ \swarrow_1, \ \nwarrow_1, \ \nwarrow_1, \ \downarrow_1, \ \uparrow_1, \ \prec_1, \ \succ_1, \ \circ_1 : E^{\otimes 2} \xrightarrow{} E$ extended to $E[[h]]$ and denoted by the tilde notation. Set:
$$ x \diamond(h) y := x \tilde{\diamond} y + hx \tilde{\diamond}_1 y ,$$
for all $\diamond \in I_E$.
Suppose the operations labelled by 1 obey $r$-ennea-algebra axioms. The sums operations can be also written, for instance $\bar{\star}(h) \longrightarrow \tilde{\bar{\star}} + h\tilde{\bar{\star}}_1,$ and so on.
We require that $(E[[h]], \ \{ \diamond(h);  \ \diamond \in I_E \})$ is a $r$-ennea-algebra. This implies that operations labelled by 1 and operations with no label have to verify a system $\widetilde{(Syst. \ 4)}$  of 49 equations, easily obtained from the definition of a $r$-ennea-algebra. These 49 equations can be solved if the involved operations come from Baxter operators and by replacing    
Condition 1 of Theorem \ref{hyp. tri3} by the following one:
\begin{enumerate}
\item {The $k$-vector space $(A, \ \mu, \ \Delta)$ and $(A, \ \mu, \ \Delta_1)$ are both $\epsilon(r)$-bialgebras.}
\end{enumerate} 
The operations we are looking for on $\textsf{End}(A)$ are defined as follows by:
$$T \searrow S := \gamma(T) \succ_{\beta} S, \ \ T \nearrow S:= T \succ_{\beta} \gamma(S), \ \ T \swarrow S := \gamma(T) \prec_{\beta} S, \ \ T \nwarrow S := T \prec_{\beta} \gamma(S),$$
$$ T \uparrow S := T \circ_{\beta} \gamma(S), \ \ T \downarrow S := \gamma(T) \circ_{\beta} S,$$
and,
$$T \searrow_1 S := \gamma(T) \succ_{\beta_1} S, \ \ T \nearrow_1 S:= T \succ_{\beta_1} \gamma(S), \ \ T \swarrow_1 S := \gamma(T) \prec_{\beta_1} S, \ \ T \nwarrow_1 S := T \prec_{\beta_1} \gamma(S).$$
$$ T \uparrow_1 S := T \circ_{\beta_1} \gamma(S), \ \ T \downarrow_1 S := \gamma(T) \circ_{\beta_1} S.$$
We observe that each equation of the system $(Syst. \ 2)$ holds, thanks to the hypothesis just above and Theorem \ref{hyp. tri2} and will give rise to 7 equations, according to the replacement of only two variables, for instance, $R$, $S$, by $\gamma(R)$, $\gamma(S)$, only one variable for instance $T$ by $\gamma(T)$ or making no replacement. This gives the expected $7 \times 7 =49$ equations.
As an example, consider the following equation from $(Syst. \ 2)$,
$$(R \circ_{\beta_1} S)\prec_{\beta} T + (R \circ_{\beta} S)\prec_{\beta_1} T = R \circ_{\beta_1}( S \prec_{\beta} T) + R \circ_{\beta}( S \prec_ {\beta_1} T),$$
for $R,S,T \in \textsf{End}(A)$.
Choose to replace only $R$ by $\gamma(R)$. This gives:
$$(\gamma(R) \circ_{\beta_1} S)\prec_{\beta} T + (\gamma(R) \circ_{\beta} S)\prec_{\beta_1} T = \gamma(R) \circ_{\beta_1}( S \prec_{\beta} T) + \gamma(R) \circ_{\beta}( S \prec_ {\beta_1} T),$$
Therefore,
$(R \downarrow_1 S)\prec_{\beta} T  + (R \downarrow S)\prec_{\beta_1} T= R \downarrow_1 (S \prec_{\beta} T) + R \downarrow (S \prec_{\beta_1} T),$
as expected. 
Choose to replace only $R$ by $\gamma(R)$ and $T$ by $\gamma(T)$. This gives:
$$(\gamma(R) \circ_{\beta_1} S)\prec_{\beta} \gamma(T) + (\gamma(R) \circ_{\beta} S)\prec_{\beta_1} \gamma(T) = \gamma(R) \circ_{\beta_1}( S \prec_{\beta} \gamma(T)) + \gamma(R) \circ_{\beta}( S \prec_ {\beta_1} \gamma(T)),$$
Therefore,
$$(R \downarrow_1 S)\nwarrow T  + (R \downarrow S)\nwarrow_1 T= R \downarrow_1 (S \nwarrow T) + R \downarrow (S \nwarrow_1 T),$$
as expected. Once this system of 49 equations is solved, we extend the operations as usual to $\textsf{End}(A)[[h]]$. A formal deformation of any enna-algebras will then be obtained provided operations comes from Baxter operators.
\Rk
A binary, quadratic and non-symmetric operad can be obtained. The sequence of these dimensions starts with $1,18,501$. By choosing the action of the unit as follows, $x \nwarrow 1 := x, 1 \searrow x := x$ or $x \nwarrow_1 1 := x, 1 \searrow_1 x := x$ (and all the other actions equal to zero) the augmented free algebra associated with this operad has a connected Hopf algebra structure. Observe also that there exists a notion of transpose and opposite associated with this type of algebras. The transformations are the same that those describing the transpose and opposite of a ennea-algebra and have to be applyied to the unlabelled operations as well as to those labelled by 1. As our choice of the unit action is in agreement with the opposite structure,
the augmented free commutative algebra will have a structure of connected Hopf algebra.
\section{Conclusion: Formal deformations of combinatorial objects}
In this paper, many binary quadratic and non-symmetric operads appear. Most of them come from formal perturbations of certain types of algebras. We conjecture that the free algebras associated with these operads can be described by combinatorial objects. As an example, it is known that the free dendriform dialgebra on one generator is described by planar binary trees. The formal deformation of dendriform dialgebras by operations verifying dendriform dialgebra axioms leads to an operad characterized by a sequence starting with 1, 4, 23. The combinatorial object describing the free algebra associated with this operad could then be viewed as a formal deformation of planar binary trees considered as dendriform dialgebra. The same thing holds for the other operads involved in this paper. What is at stake is a classification and a deformation of known combinatorial objects viewed as free algebras of certain types of algebras.

\noindent
\textbf{Acknowledgments:}
The author is indebted to J.-L. Loday for giving him articles prior to their
publications. He thanks M. Aguiar and J.-L. Loday for having read preliminary drafts of this paper.

\bibliographystyle{plain}
\bibliography{These}

\begin{thebibliography}{10}

\bibitem{Aguiar}
M.~{\textsc{Aguiar}}.
\newblock Infinitesimal bialgebras, pre-lie and dendriform algebras.
\newblock {\em to appear in ``{H}opf algebras: {P}roceedings from an
  {I}nternational {C}onference held at {D}e{P}aul {U}niversity'';
  arXiv:math.QA/0211074}.

\bibitem{AguiarLoday}
M.~{\textsc{Aguiar}} and J.-L. {\textsc{Loday}}.
\newblock Quadri-algebras.
\newblock {\em arXiv:math.QA/0309171}.

\bibitem{Baxter}
G.~{\textsc{Baxter}}.
\newblock An analytic problem whose solution follows from a simple algebraic
  identity.
\newblock {\em Pacific J. Math.}, 10:731--742, 1960.

\bibitem{CKreimer}
A.~{\textsc{Connes}} and D.~{\textsc{Kreimer}}.
\newblock Renormalisation in quantum field theory and the {R}iemann-{H}ilbert
  problem.
\newblock {\em eprint arXiv:hep-th/9909126}.

\bibitem{Domb}
C.~{\textsc{Domb}} and A.~J. {\textsc{Barrett}}.
\newblock Enumeration of ladder graphs.
\newblock {\em Discrete Math.}, pages 341--358, 1974.

\bibitem{KEF}
K.~{\textsc{Ebrahimi-Fard}}.
\newblock Loday-type algebras and the {R}ota-{B}axter relation.
\newblock {\em eprint, math-ph/0207043}.

\bibitem{Flajolet}
P.~{\textsc{Flajolet}} and M.~{\textsc{Noy}}.
\newblock Analytic combinatorics of non-crossing configurations.
\newblock {\em Discrete Math.}, 204:203--229, 1999.

\bibitem{Fresse}
B.~{\textsc{Fresse}}.
\newblock Koszul duality of operads and homology of partition posets.
\newblock {\em Preprint 2002}.

\bibitem{Gerst}
M.~{\textsc{Gerstenhaber}}.
\newblock The cohomology structure of an associative ring.
\newblock {\em Ann. of Math.}, 78:267--288, 1963.

\bibitem{Gerst1}
M.~{\textsc{Gerstenhaber}} and J.D.~{\textsc{Stasheff}} (Eds.).
\newblock {\em Deformation theory and quantum groups with applications to
  mathematical physics}, volume 134.
\newblock Proceedings of an AMS-IMS-SIAM 1990 joint summer research conference,
  Contemporary Mathematics, 1992.

\bibitem{GK}
V.~{\textsc{Ginzburg}} and M.~{\textsc{Kapranov}}.
\newblock Koszul duality for operads.
\newblock {\em Duke Math. J.}, 76(1):203--272, 1994.

\bibitem{Rota}
S.A. {\textsc{Joni}} and G.-C. {\textsc{Rota}}.
\newblock Coalgebras and bialgebras in combinatorics.
\newblock {\em Stud. Appl. Math.}, 61:93--139, 1979.

\bibitem{Kreimer}
D.~{\textsc{Kreimer}}.
\newblock Chen's iterated integral represents the operator product expansion.
\newblock {\em eprint arXiv:hep-th/9901099, Adv. Theor. Math. Phys},
  3:627--670, 1999.

\bibitem{BaxLer}
Ph. {\textsc{Leroux}}.
\newblock On some remarkable operads constructed from {B}axter operators.
\newblock {\em In preparation}.

\bibitem{codialg1}
Ph. {\textsc{Leroux}}.
\newblock Tiling the $(n^2,1)$-{D}e-{B}ruijn graph with $n$ coassociative
  coalgebras.
\newblock {\em eprint arXiv:math.QA/ 0209108}.

\bibitem{Lodayscd}
J.-L. {\textsc{Loday}}.
\newblock Scindement d'associativit\'e et alg\`ebres de {H}opf.

\bibitem{LodayLeib}
J.-L. {\textsc{Loday}}.
\newblock Une version non commutative des alg\`ebres de {L}ie: Les alg\`ebres
  de {L}eibniz.
\newblock {\em L'Enseignement Math.}, 39:269--293, 1993.

\bibitem{Loday}
J.-L. {\textsc{Loday}}.
\newblock Dialgebras.
\newblock {\em in Dialgebras and related operads, Lecture Notes in Math.},
  1763:7--66, 2001.

\bibitem{LodayRonco}
J.-L. {\textsc{Loday}} and M.~{\textsc{Ronco}}.
\newblock Trialgebras and families of polytopes.
\newblock {\em eprint arXiv:math.QA/0205043}.

\bibitem{Richter}
B.~{\textsc{Richter}}.
\newblock Dialgebren, {D}oppelagebren und ihre {H}omologie.
\newblock {\em {D}iplomarbeit, {B}onn {U}niversit\"at, unpublished}.

\bibitem{Rota1}
G.-C. {\textsc{Rota}}.
\newblock Baxter operators, an introduction.
\newblock {\em in {G}ian-{C}arlo {R}ota on {C}ombinatorics: {I}ntroductory
  papers and commentaries ({J}.{P}.{S}. {K}ung, {E}d.) {B}irkhauser, {B}oston},
  pages 504--512, 1995.

\bibitem{Voi}
D.~{\textsc{Voiculescu}}.
\newblock Free analysis question {I}: Duality transform for the coalgebra of
  $\partial_{X:B}$.
\newblock {\em arXiv:math.QA/0306172}.

\end{thebibliography}

\end{document}